\documentclass[11pt]{article}

 \usepackage{mst-stylefile}
 \usepackage{Bbb11a}
 \usepackage{harvard}
 \usepackage{amsfonts}
 \usepackage{imsart}

 \usepackage{graphicx}

 \usepackage{color}

\include{epsf}

\renewcommand\cal{\mathcal }

\renewcommand{\Bbb}{\mathbb}

\newcommand{\bbR}{{\Bbb R}}

\newcommand{\bbN}{{\Bbb N}}

\newcommand\const{\mbox{\rm const}\,}

\newcommand\E{\Bbb E}
\renewcommand\P{\Bbb P}
\newcommand\noi{\noindent}

\newcommand{\what}{\widehat}

\newcommand{\cov}{\mbox{{\rm Cov}}}

\def\D{{\cal D}([0,\infty),\bbR)}

\def\Sum{\mathop{\sum}}

\def\wtilde{\widetilde }

\renewcommand{\cite}{\citeyear}

\newtheorem{conjecture}{\bf Conjecture}[section]

\def\bmu{\mu}

\def\btheta{\theta}
\def\bY{Y}

\def\D{D} 

\begin{document}
\begin{frontmatter}

\title{Estimating heavy--tail exponents through max self--similarity}
\runtitle{Estimating heavy--tail exponents through max self--similarity}
\author{Stilian A.\ Stoev\thanks{{\it Corresponding author:} {\tt sstoev@umich.edu},
  Department of Statistics, University of Michigan,
 439 West Hall, 1085 South University,
 Ann Arbor, MI 48109-1107, U.S.A. phone: (734) 763-3519, fax: (734) 763-4676}}
\address{University of Michigan, Ann Arbor}
\author{George Michailidis}
\address{University of Michigan, Ann Arbor}
\author{Murad S.\ Taqqu}
 \address{Boston University}
\runauthor{Stoev {\it et  al.\ }}

\begin{abstract}
In this paper, a novel approach to the problem of estimating the heavy--tail exponent 
$\alpha>0$ of a distribution is proposed.  It is based on the fact that block--maxima of size $m$
of the independent and identically distributed data scale at a rate of $m^{1/\alpha}$.
This scaling rate can be captured well by the {\it max--spectrum} plot of the data that leads
to regression based estimators.
Consistency and asymptotic normality of these estimators is 
established under mild conditions on the behavior of the tail of the distribution. The 
results are obtained by establishing bounds on the
rate of convergence of moment--type functionals of heavy--tailed maxima.  Such
bounds often yield exact rates of convergence and are of independent interest.
Practical issues on the automatic selection of tuning parameters for the estimators 
and corresponding confidence intervals are also addressed.  Extensive numerical
simulations show that the proposed method proves competitive for both small and
large sample sizes and for a large range of tail exponents. The method is shown to be more 
robust than the classical Hill plot and is illustrated on two data sets of
insurance claims and natural gas field sizes.
\end{abstract}



\end{frontmatter}
\doublespacing
\section{Introduction}
 \label{s:intro}
 
Heavy--tailed distributions arise in many diverse scientific areas:
insurance claims, high--speed network traffic, hydrology, the topological structure of 
the World Wide Web and of social networks, linguistics, just to name a few (see e.g.\ 
Adler et al.\ \cite{adler:feldman:taqqu:1998}, McNeil \cite{mcneil:1997}, Resnick \cite{resnick:1997}, 
 Faloutsos et al.\ \cite{faloutsos:faloutsos:faloutsos:1999}, Adamic and
 Huberman \cite{adamic:huberman:2000,adamic:huberman:2002},
 Zipf \cite{zipf:1932,zipf:1949}, Tsonis et al.\
\cite{tsonis:schultz:tsonis:1997}).  Highly optimized physical systems also exhibit heavy--tailed behavior, 
as discussed in Carlson and Doyle \cite{carlson:doyle:1999}.

A real valued random variable $X$ with cumulative distribution
function (c.d.f.)\ $F(x) = \P\{X\le x\},\ x\in\bbR$ is said to have
(right) heavy tail if, 
\beq\label{e:F}
 \P\{ X>x \} = 1-F(x) = L(x) x^{-\alpha}, \ \mbox{ as }x\to\infty
\eeq
for some $\alpha>0$, where $L(x)>0$ is a slowly varying function.  The {\em tail exponent}
$\alpha>0$ controls the rate of decay of $F$ and hence characterizes
its tail behavior. The problem of estimating the tail exponent has
attracted a lot of attention in the literature since it poses
numerous theoretical, as well as, practical challenges (de
Haan et al. \cite{dehaan:drees:resnick:2000} and de Sousa
and Michailidis \cite{desousa:michailidis:2004}). Most
approaches focus on the scaling behavior of the largest order
statistics $X(1;N) \ge X(2;N) \ge \cdots \ge X(N;N)$ obtained from
an independent and identically distributed (i.i.d.)\ sample $X(1),\ldots,X(N)$ from $F$.  
Typical examples include Hill's estimator \cite{hill:1975}, its numerous variations (Kratz and Resnick \cite{kratz:resnick:1996},
Resnick and St{\v{a}}ric{\v{a}} \cite{resnick:starica:1997}), and the kernel--based estimators of Cs\"org{\H{o}} et al.
\cite{csorgo:deheuvels:mason:1985} (see also Feuerverger and Hall \cite{feuerverger:hall:1999}).  For example, 
the Hill estimator, which is one of the most widely used estimators in practice, can be written as
\beq\label{e:alpha-k-Y}
 \what\alpha_H(k) = {\Big(}\frac{1}{k} \sum_{i=1}^k i (\ln X(i;N) - \ln X(i+1;N))
  {\Big)}^{-1}
   =: {\Big(} \frac{1}{k} \sum_{i=1}^k Y_i {\Big)}^{-1},
\eeq
where $Y_i:= i (\ln X(i;N) - \ln X(i+1;N))$. As shown in Weissman \cite{weissman:1978}, 
assumption \refeq{F} implies that for all fixed $k$'s, the vector
$\{Y_i\}_{i=1}^k$ converges in distribution to a vector of independent exponentially 
distributed variables with mean $1/\alpha$.  Therefore, when both $N$ and $k$ are large,
the statistic $\what\alpha_H(k)$ in \refeq{alpha-k-Y} behaves like the sample mean of 
a sample of independent exponential variables. This suggests that the estimator $\what\alpha_H(k)$
is consistent (Mason \cite{mason:1982}), and under some additional conditions on the tail
behavior of $F$, asymptotically normal (Hall \cite{hall:1982}). 
 In practice, one relies on plotting $\what\alpha_H(k)$ as a function of
the order statistics $k$ (Hill plot) and then selecting an appropriate value for
$k$ (see example in Figure \ref{fig:first-data-example}). In the case of the Pareto distribution
($F(x) = 1- (x/\sigma_0)^{-\alpha},\ x\ge \sigma_0,~\sigma_0>0$), the
Hill estimator is also a conditional maximum likelihood estimator.
However, when deviations from this ideal case occur, it exhibits substantial bias and the
resulting plot can be misleading (see examples and discussion in 
\ de Haan et al. \cite{dehaan:drees:resnick:2000} and de Sousa and Michailidis
\cite{desousa:michailidis:2004} and references therein). These shortcomings were addressed
in a series of papers that introduced modifications of the original Hill estimator and the
resulting Hill plot. The kernel--type estimators introduced by 
Cs\"org{\H{o}} et al.\ \cite{csorgo:deheuvels:mason:1985} 
extend the Hill estimator, by introducing non--uniform weights in \refeq{alpha-k-Y} (see
also Groeneboom {\it et al.\ }\cite{groeneboom:lopuhaa:dewolf:2003}).  
Namely, given a non--negative and non--increasing kernel function $K(x),\ x>0$, one considers
\beq\label{e:kernel}
 \what \alpha_{K,\lambda, N} := {\Big(}\frac{1}{N}\sum_{i=1}^N K(i/\lambda N) Y_i {\Big)}^{-1} \int_0^{1/\lambda} K(x) dx,
\eeq
for some $\lambda>0$.  The Hill estimator can be recovered as a special choice of 
the function $K$. Observe also that the threshold parameter $k$ in \refeq{alpha-k-Y} is 
no longer present.  The choice of the kernel function and the bandwidth
parameter $\lambda>0$, however, remain an important and difficult problem for the
kernel estimators, similar to the choice of $k$ for the Hill estimator.  One practical
disadvantage of kernel--type estimators is that no analogue of the Hill plot exists.
Therefore, one cannot readily judge how reliable the resulting numerical estimates are. 

Other important and popular estimators include the Pickands estimator  (see, Pickands \cite{pickands:1975} and 
Dekkers and de Haan \cite{dekkers:dehaan:1989}) and de Haan's moment type estimator
(see Dekkers et al. \cite{dekkers:einmahl:dehaan:1989}). 
Resnick and St{\v{a}}ric{\v{a}} \cite{resnick:starica:1997} introduced a modified and smoothed version
of the Hill plot and showed that it performs better in practice when the data depart from the Pareto model
(see also de Haan et al. \cite{dehaan:drees:resnick:2000}).  The consistency of estimators based
on this alternative Hill plot is also established for dependent data (see, Resnick and
St{\v{a}}ric{\v{a}} \cite{resnick:starica:1995}).

\medskip
In this study, we propose a novel method for estimating the tail
index $\alpha$. It relies on the concept of {\em max
self-similarity}.  We focus on the case when the slowly
varying function in \refeq{F} is asymptotically constant and consider
block--wise maxima of i.i.d.\ random variables $X(1),X(2),\ldots$ with
c.d.f.\ $F$.  Block--maxima of block sizes $m$, scale at a rate of
$m^{1/\alpha}$, as $m\to\infty$.  Therefore, we can obtain an estimate
of $\alpha$, by focusing on a sequence of growing, dyadic block sizes
$m = 2^j,\ 1\le j \le \log_2 N,\ j\in\bbN$, and estimating the mean of
logarithms of block--maxima (log--block--maxima).  This is achieved by
examining the {\it max--spectrum plot} of the data, defined as means of
log--block--maxima as a function of the logarithm of the
block--size. The slope of the max--spectrum plot for large
block--sizes yields an estimate of $1/\alpha$ (see Figure
\ref{fig:first-data-example} below).

When the $X(i)$'s come from a Fr\'echet distribution, then their block--maxima have the
same Fr\'echet distribution, rescaled by $m^{1/\alpha}$, where $m$ denotes the block size.   Thus, 
in practice, the max--spectrum plot is essentially linear (Figure \ref{fig:max_spectrum}). One can view 
i.i.d.\ Fr\'echet sequences
as {\it max self--similar} with self--similarity parameter $1/\alpha$ (Definition \ref{d:max-ss}).
Due to this exact max self--similarity property, our estimation framework works best for Fr\'echet data.  
On the other hand, the Hill--type estimators work best for Pareto data.  
This also shows the fundamental difference between the two approaches.
In many important applications the Hill plot is rather volatile.  The max spectrum turns out to be
more robust to outliers in the data or to deviations from its corresponding ideal Fr\'echet model than the 
Hill plot.  In Section \ref{s:data}, we examine two data sets: {\it (i)} $2,167$ insurance claims 
due to fire losses in Denmark and {\it (ii)} volumes of natural gas reserves in $406$ Oil rich provinces.
In both cases, the max self--similarity estimators yield values consistent with previous detailed
studies of these data sets (see McNeil \cite{mcneil:1997} and de Sousa and Michailidis
\cite{desousa:michailidis:2004}, respectively).  These values depart from values that one obtains
directly from the Hill plots.  In fact, in case {\it (ii)}, due to the peculiar discrete nature 
of the data set the Hill plot has a saw tooth shape and it is particularly hard to interpret, whereas the
max spectrum plot appears to yield a reliable estimate.

\medskip
The remainder of the paper is structured as follows.
In Section \ref{s:max-ss}, we introduce the max--spectrum plot and the
self--similarity estimators of the heavy--tail exponent $\alpha$ and establish their
basic properties in the ideal Fr\'echet setting. Some useful results on rates for
moment--type functionals of heavy--tailed maxima are presented in Section \ref{s:rate}.  These results
are used to prove the consistency and asymptotic normality of the max self--similarity estimators 
in Section \ref{s:cons-an}. In Section \ref{s:performance}, the performance of
the new estimators is examined through a simulation study.  The max self--similarity
estimators are then shown to work well in the context of two challenging real data examples 
where the classical Hill plot is rather volatile and is hard to interpret.

\section{Max self--similarity and tail exponent estimators}
 \label{s:max-ss}

 In this section, we introduce some notation and recall some basic
 definitions used in the remainder of the paper.  We then introduce estimators
 of the heavy--tail exponents based on max self--similarity and discuss 
 their basic properties in the ideal Fr\'echet  case.

 \subsection{Definition and basic properties}
 
We focus on the case 
where the slowly varying function $L$ in \refeq{F} is trivial, that is, when
\beq\label{e:F-tail}
  \P\{X>x\} = 1 - F(x) \sim \sigma_0^\alpha x^{-\alpha},\ \mbox{ as } x\to\infty, 
\eeq
with $\sigma_0>0$ and where $\sim$ means that the ratio of the left--hand side (l.h.s.)
to the right--hand side (r.h.s.) in \refeq{F-tail} tends to $1$, as $x\to\infty$.
For simplicity, we further assume that the $X(i)$'s are almost surely positive ($F(0)=0$).
We address the general case where the $X(i)$'s can take negative values in Section \ref{s:cons-an} (see, 
Proposition \ref{p:coupling}).

We begin with some useful definitions: for an i.i.d.\ sample $X(i),\ i\in\bbN:=\{1,2,\ldots\}$ from $F$,
consider the sequence of block--maxima
$$
   X_m(k) := \max_{1\le i \le m}X(m(k-1)+i) \equiv \bigvee_{i=1}^m X(m(k-1)+i),\ \ k=1,2,\ldots,
$$
with $m\in\bbN$, where $X_m(k)$ is the greatest observation in the $k-$th block. The Fisher--Tippett--Gnedenko Theorem
(see e.g.\ Proposition 0.3 in Resnick \cite{resnick:1987}) then implies that, as $m\to\infty$, $m^{-1/\alpha} X_m(k)$
converges in distribution to a random variable $Z$ with an $\alpha-$Fr\'echet
distribution.  More precisely,
\beq\label{e:Z}
 \P\{Z\le x\} = \exp\{-\sigma_0^\alpha x^{-\alpha}\},\ \ x>0,
\eeq
where $\sigma_0>0$, called the {\it scale coefficient} of $Z$, is as in \refeq{F-tail}.  In fact, as
$m\to\infty$, we have
\beq\label{e:X_m}
   {\Big\{}\frac{1}{m^{1/\alpha}} X_m(k){\Big\}}_{k\in\bbN} \stackrel{d}{\longrightarrow}
   {\Big\{} Z(k) {\Big\}}_{k\in\bbN},
\eeq
where the $Z(k)$'s are independent copies of $Z$ and where $\stackrel{d}{\to}$ denotes convergence of 
the finite--dimensional distributions.  Thus, for large values of $m$, the normalized block--maxima 
behave like a sequence of i.i.d.\  $\alpha-$Fr\'echet variables.  In fact, when the $X(k)$'s are $\alpha-$Fr\'echet,
\refeq{X_m} holds with equality for all $m\in\bbN$ (see Relation \refeq{Z-max-stable} in the Appendix).
The sequence of i.i.d.\  $\alpha-$Fr\'echet $X(k)$'s is thus {\it max self--similar} in the sense 
of the following definition.

\begin{definition}\label{d:max-ss}
 A sequence of random variables $X = \{X(k)\}_{k\in\bbN}$ (defined on the same probability space)
 is said to be max self--similar with self--similarity parameter $H>0$, if for any $m>0,\ m\in\bbN$,
\beq\label{e:d:max-ss}
 {\Big\{}\bigvee_{i=1}^{m} X(m(k-1)+i) {\Big\}}_{k\in\bbN} \stackrel {d}{=} 
 {\Big\{}m^H X(k){\Big\}}_{k\in\bbN},
\eeq
where $=^d$ denotes equality of the finite--dimensional distributions.
\end{definition}

If the $X(k)$'s are i.i.d.\  but not Fr\'echet, then 
Relation \refeq{X_m} indicates that \refeq{d:max-ss} holds asymptotically, as $m\to\infty$, with $H=1/\alpha$. 
Thus, any sequence of i.i.d.\  heavy--tailed variables can be regarded as {\it asymptotically max
self--similar} with self--similarity parameter $H=1/\alpha$.  This feature suggests that an estimator of $H$ and 
therefore $\alpha$ can be obtained by focusing on the scaling of the block--maxima of growing block sizes.
Crovella and Taqqu \cite{crovella:taqqu:1999} used a similar idea based on the scaling of block--wise sums
to estimate a heavy--tail exponent $\alpha$ when $\alpha \in (0,2)$.

\medskip
Given an i.i.d.\  sample $X(1),\ldots,X(N)$ from $F$, we consider 
\beq\label{e:C-j}
  \D(j,k) := \max_{1\le i \le 2^j} X(2^j(k-1)+i) = \bigvee_{i=1}^{2^j} X(2^j(k-1)+i),\ \ k = 1,2,\ldots, N_j,
\eeq 
for all $j=1,2,\ldots,[\log_2 N],$  where $N_j := [N/2^j]$ and $[x]$ denotes the largest integer not
greater than $x\in\bbR$.  
By analogy to the discrete wavelet transform, we refer to the parameter $j$ as the {\it scale} and
to $k$ as the {\it location} parameter.  We consider dyadic block--sizes for
algorithmic and computational convenience (for more details, see Stoev et al.\ \cite{stoev:michailidis:hamidieh:taqqu:2006P}).

Observe that for any fixed $j$, the block--maxima $\D(j,k)$ are independent in $k$ since they involve maxima over 
non--overlapping blocks of the $X(i)$'s.  Moreover, as argued above, when the $X(i)$'s follow
an $\alpha-$Fr\'echet distribution,
\beq\label{e:Cj}
  \{\D(j,k)\}_{k\in\bbN} \stackrel{d}{=} \{2^{j/\alpha} \D(0,k)\}_{k\in\bbN} = \{2^{j/\alpha} X(k)\}_{k\in\bbN},
\eeq
for any scale $j \in\bbN$.  Introduce the statistics
\beq\label{e:Yj}
  Y_j := \frac{1}{N_j} \sum_{k=1}^{N_j} \log_2 \D(j,k),\ \ j=1,2,\ldots,[\log_2(N)]
\eeq
and observe that by the Law of Large Numbers, the $Y_j$'s are consistent and unbiased estimators of
the expectations $\E \log_2 \D(j,1)$, provided that these are finite.  (Corollary \ref{c:logs}
below establishes that $\E|\log_2 \D(j,1)|$ are finite under general conditions
on the c.d.f.\  $F(x)$.)  
In view of the asymptotic max self--similarity \refeq{X_m} of $X$, relationship \refeq{Cj} holds approximately
for large scales $j$, and in fact,
\beq\label{e:Yj-sim}
  \E Y_j = \E \log_2 \D(j,1) \simeq j/\alpha + C, 
\eeq
with $C = C(\sigma_0,\alpha) = \E \log_2 \sigma_0 Z$, where $Z$ is an $\alpha-$Fr\'echet variable with unit coefficient
as in \refeq{Z} above. Here $\simeq$ means that the difference between the l.h.s.\ and the r.h.s.\ tends to zero. 

In practice, one can look at the {\it max--spectrum plot} of the statistics $Y_j$'s versus $j$ 
(see Figure \ref{fig:first-data-example} below).  In view of \refeq{Yj-sim} it is expected that for
large $j$'s the slope coefficient of a linear fit of the $Y_j$'s  against $j$'s would yield an estimate of $H = 1/\alpha$.   
Further, observe that the log--linear scaling relation in \refeq{Yj-sim} becomes more 
precise, the larger the scale $j$ (block--size $2^j)$ and holds exactly for all scales
$j=1,\ldots,[\log_2(N)]$, when the $X(k)$'s come from an $\alpha-$Fr\'echet distribution (see \refeq{Cj}).

\medskip
Thus, given a range of scales $1\le j_1\le j \le j_2 \le [\log_2(N)]$, we define the following
regression--based estimators of $H = 1/\alpha$ and $\alpha$
\beq\label{e:H-hat}
 \what H_w(j_1,j_2) := \sum_{j=j_1}^{j_2} w_j Y_j,\ \ \mbox{ and }\ \ \what \alpha_w (j_1,j_2) := 1/\what H_w(j_1,j_2),
\eeq
where the weights $w_j$ are chosen so that
\beq\label{e:wj}
 \sum_{j=j_1}^{j_2} w_j = 0\ \ \ \mbox{ and } \ \ \ \sum_{j=j_1}^{j_2} j w_j = 1.
\eeq
It is easy to see that the linear estimators $\what H_w$ in \refeq{H-hat} with weights as in \refeq{wj} are
least squares estimators in a linear regression model.  In the rest of the paper, the estimators 
$\what H_w$ and $\what \alpha_w$ in \refeq{H-hat} are referred to as {\it max self--similarity} estimators.

\medskip
\noindent
{\bf Remark} {\it (Computational complexity)} \\
The proposed estimators exhibit a significant computational advantage over 
Hill--type or kernel--based estimators.  Given a sample of size $N$ one can
compute the max--spectrum $Y_j,\ 1\le j \le [\log_2 N],$ with $Y_j$ as
in \refeq{Yj} by using ${\cal O}(N)$ operations since 
${\cal O}(N/2^j)$ pair--wise maxima and sums are computed, for $j=1,\ldots,[\log_2
N]$, and therefore
$
{\cal O} {\Big(} \sum_{j=1}^{[\log_2 N]} [N/2^j] {\Big)} = {\cal O}(N)
$ 
operations are done. On the other hand, 
methods involving order statistics require sorting the
sample which results in ${\cal O}(N \log_2(N))$ operations.
\medskip

We now illustrate the nature of the max-spectrum plot and the resulting estimator
using an example of Internet topology data. The data
describe the degree of connectivity between autonomous systems (AS -
networks under a single administrative authority) on the Internet for
the year 2002 and is provided by the National Laboratory for Applied
Network Research. The information has been used to characterize the
topology of the Internet (see, e.g.\ Faloutsos et al.\
\cite{faloutsos:faloutsos:faloutsos:1999} and Chen et al.\
\cite{chen:chang:govindan:jamin:shenker:willinger:2002}).  The size of
the data set is 13,579 and each observation gives the number of
connections of an AS to peer AS.  The histogram of the data (in
log-scale) shows that the vast majority of the AS are connected to very
few peer systems, but there are a few AS that are directly connected to over
10\% of their peer systems. The max--spectrum indicates a value for the
tail index of about 1.5. The Hill estimator for $k=80$ (where the Hill
plot seems to stabilize) suggests a value of 1.43.

\begin{figure}[h!]
\begin{center}
\includegraphics[height=2in,width=2.5in]{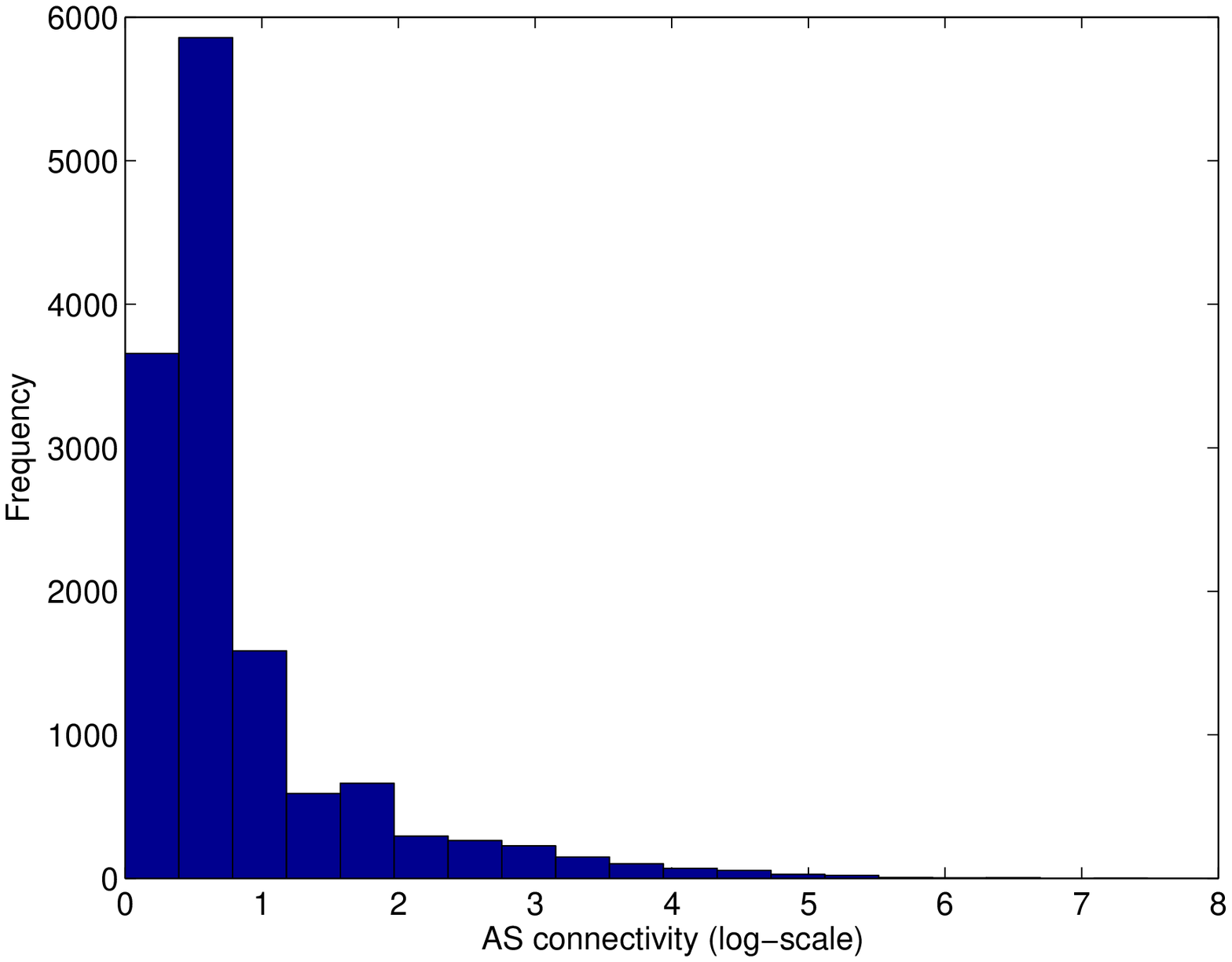} \hskip .2 in
\includegraphics[height=2in,width=2.5in]{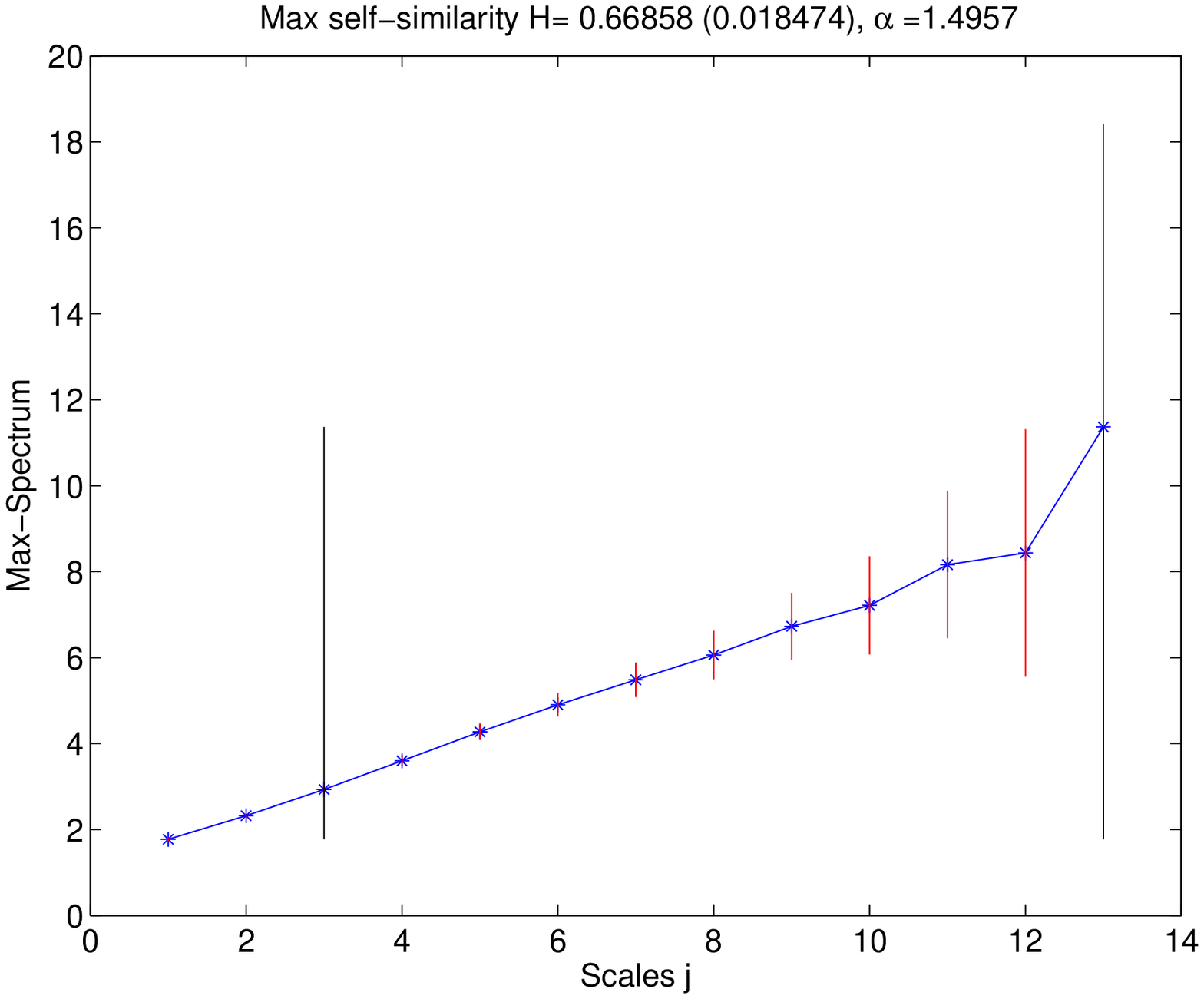}
{\caption{\label{fig:first-data-example} \small {\it Left panel:} histogram (log--scale) of AS
connectivities. {\it Right panel:} max--spectrum plot for the AS connectivity data. 
The large vertical lines indicate the range of $j$'s where a linear fit was  
used to estimate the heavy--tail index $\alpha$.  The shorter vertical lines are $95\%$ 
confidence intervals for the $\E Y_j$'s.
The reciprocal of the slope yields an estimate of $\what \alpha_w(3,13) =  
1.4957$.
This range was selected automatically with tunning level $p=0.1$, discussed in Section \ref{s:cut-off}.
}}
\end{center}
\end{figure}

\subsection{The ideal Fr\'echet case}
 \label{s:GLS}

We start by assuming that $X(1),\ldots,X(N)$ is an i.i.d.\  sample
of $\alpha-$Fr\'echet variables with scale coefficient $\sigma_0>0$ and study the behavior of 
$\what H_w(j_1,j_2)$ in this setting.

\medskip
Consider the regression problem
\beq\label{e:Yj-epsilonj}
Y_j =  j/\alpha + C + \epsilon_j,\ \ j_1\le j \le j_2
\eeq
where 
\beq\label{e:the-C}
  C = C(\sigma_0,\alpha) = \E \log_2(\sigma_0Z) = \log_2(\sigma_0) + \E\log_2(Z)
\eeq
for an $\alpha-$Fr\'echet $Z$ random variable 
with unit scale coefficient, and where $1\le j_1\le j_2\le [\log_2
N]$. In view of \refeq{Cj}, we have that the errors $\epsilon_j$ have
zero means.  They are, however, dependent in $j$ due to the
corresponding dependence of the $Y_j$ statistics in \refeq{Yj}.
Moreover, the number of $\D(j,k)$'s at a scale $j$ in \refeq{Yj} is
$N_j= [N/2^j]$ and therefore, the variances of the $\epsilon_j$'s grow
exponentially in $j$.  This implies that the minimal variance unbiased estimators 
of the parameters of interest $\theta = (H,C)^t$ that are linear in
$Y_j$ are obtained through {\it generalized least squares} (GLS). They are given by
\beq\label{e:theta-GLS}
 \what \theta_{\Sigma} = {\what H_{\Sigma} \choose \what C_{\Sigma} }
  =(A^t \Sigma^{-1} A)^{-1} A^t \Sigma^{-1} Y,
\eeq
where $A = (a\, b)$ with $a^t = (j_1,\ldots, j_2)$ and $b^t = (1,\ldots, 1)$, and
$\Sigma = (\cov(Y_i,Y_j))_{i,j=j_1}^{j_2}$  is the covariance matrix of the vector 
$Y=\{Y_j\}_{j=j_1}^{j_2}$. 
An explicit expression of the matrix $\Sigma= \Sigma_\alpha(j_1,j_2;N)$ is given next.

\begin{proposition}
\label{p:cov}  Let $Y=\{Y_j\}_{j=j_1}^{j_2}$ be as in \refeq{Yj}, where the underlying distribution of the
 $X(k)$'s is $\alpha-$Fr\'echet with scale coefficient $\sigma_0>0$.  Then, for all $j_1\le i\le j\le j_2$,
$$
 \E Y_j = j/\alpha + C(\sigma_0,\alpha),
$$
and 
\beq\label{e:p:cov}
  \cov(Y_i,Y_j) = \Sigma_\alpha(j_1,j_2;N)_{ij} =  \frac{2^{j-i}}{\alpha^2 N_i} \psi(|i-j|), \ \ N_i = [N/2^i],
\eeq
where 
\beq\label{e:psi}
 \psi(a) := \cov(\log_2(Z_1), \log_2(Z_1 \vee (2^a -1)Z_2)),\ a\ge 0,
\eeq 
and where $Z_1$ and $Z_2$ are independent $1-$Fr\'echet variables with 
unit scale coefficients.
\end{proposition}
\begin{proof}  Let $j_1\le i<j\le j_2$ and observe that $N_i = 2^{j-i}N_j + R$, where $0\le R<2^{j-i},\ R\in\bbN$.
 In view of \refeq{Yj},
 \begin{eqnarray}\label{e:p:cov-0.5}
  \cov(Y_i,Y_j)& = &\frac{1}{N_i N_j} \Sum_{k_1=1}^{N_i} \sum_{k_2=1}^{N_j} \cov(\log_2 \D(i,k_1), \log_2 \D(j,k_2)) \nonumber\\
               & = & \frac{1}{N_i N_j} \Sum_{k_1 =1}^{N_j}\sum_{\ell=1}^{2^{j-i}} \sum_{k_2=1}^{N_j}
 \cov(\log_2 \D(i,(k_1-1)2^{j-i}+\ell), \log_2 \D(j,k_2))\nonumber\\
 & & \ \ \ \ \ \  + \frac{1}{N_i N_j} \Sum_{\ell=1}^{R} \sum_{k_2=1}^{N_j} \cov(\log_2 \D(i,N_j2^{j-i}+\ell), \log_2 \D(j,k_2)),
 \end{eqnarray}
where the last relation follows from expressing the sum $\Sum_{k_1=1}^{N_i}$ as a double sum $\Sum_{k_1 =1}^{N_j}\sum_{\ell=1}^{2^{j-i}}$ plus
the remainder term $\Sum_{\ell=1}^{R}\sum_{k_2=1}^{N_j}$.  Observe that in view of \refeq{C-j}, we have that the terms
$\cov( \log_2 \D(i,(k_1-1)2^{j-i}+\ell) , \log_2 \D(j,k_2)),\ 1\le \ell < 2^{j-i}$ are non--zero only if $k_1 = k_2$
since otherwise the terms $\D(i,(k_1-1)2^{j-i}+\ell)$ and $\log_2 \D(j,k_2)$ involve maxima of non--overlapping sets of
$X(k)$'s. Note moreover that
\beq\label{e:p:cov-1}
 \D(j,k_2) = \D(i,(k_2-1)2^{j-i} + 1) \vee \cdots \vee \D(i,k_2 2^{j-i}),
\eeq
where the $\D(i,k)$'s are i.i.d.\  $\alpha-$Fr\'echet variables with scale coefficient $2^{i/\alpha}\sigma_0$ (see \refeq{Z-max-stable} 
below).  Therefore, for all $k=1,\ldots,N_{j}$ and $\ell = 1,\ldots,2^{j-i}$,
$$
 (\D(i,(k-1)2^{j-i}+\ell), \D(j,k)) \stackrel{d}{=} (2^{i/\alpha}Z', 2^{i/\alpha}Z' \vee (2^{j/\alpha} - 2^{i/\alpha})Z''),
$$
where $Z'$ and $Z''$ are independent $\alpha-$Fr\'echet variables with scale coefficients $\sigma_0>0$.  
Observe that $Z' = \sigma_0 Z_1^{1/\alpha}$, where $Z_1$ is $1-$Fr\'echet with unit scale coefficient. Hence, for all 
$k_1=k_2 =1,\ldots,N_j$ and $\ell =1,\ldots,2^{j-i}$, we have
\begin{eqnarray}\label{e:p:cov-2}
& & \cov( \log_2 \D(i,(k_1-1)2^{j-i}+\ell) , \log_2 \D(j,k_2)) \nonumber\\
 & &\ \ \ \ \  \ \ \ = \cov{\Big(}\log_2(2^{i/\alpha}\sigma_0 Z_1^{1/\alpha}), 
\log_2( 2^{i/\alpha}\sigma_0 Z_1^{1/\alpha} \vee (2^{j/\alpha}-2^{i/\alpha})\sigma_0 Z_2){\Big)}\nonumber\\
 &  &\ \ \ \ \ \ \ \ = \cov {\Big(}\log_2(Z_1^{1/\alpha}), \log_2(Z_1^{1/\alpha}\vee (2^{(j-i)/\alpha}-1) Z_2^{1/\alpha}){\Big)}
 =\frac{1}{\alpha^{2}} \psi(|i-j|).
\end{eqnarray}
The last two relations follow from the facts that
$\log_2(2^{i/\alpha}\sigma_0 Z_1^{1/\alpha})$ equals $\log_2(2^{i/\alpha}\sigma_0) + \alpha^{-1}\log_2(Z_1)$ and
since $\cov(\xi + a, \eta +b)  = \cov(\xi,\eta)$, for any constants $a$ and $b$ and random variables $\xi$ and $\eta$
with finite variance.

Note that the covariances in the remainder term in \refeq{p:cov-0.5} vanish since $\D(i,N_j 2^{j-i}+\ell),\ \ell=1,\ldots,2^{j-i}$
are independent of $X(i),\ i=1,\ldots, N_j2^{j}$.  Thus, by using Relation \refeq{p:cov-2}, we obtain \refeq{p:cov}.  
$\Box$
\end{proof}

\medskip
\noi{\bf Remarks}
\begin{enumerate}
 \item  Observe that the covariance matrix $\Sigma$ does not
depend on the scale coefficient $\sigma_0$, which is due to
the fact that the $Y_j$'s are obtained through a logarithmic transformation of the $X(k)$'s.

 \item Observe that for all $1\le j_1<j_2\le [\log_2 N]$ and $\alpha>0$, we have by \refeq{p:cov} that
$$ 
  \Sigma_\alpha(j_1,j_2;N) = \frac{1}{\alpha^2} \Sigma_1(j_1,j_2;N),
$$
where $\Sigma_1(j_1,j_2;N)$ corresponds to the covariance matrix of $Y=\{Y_j\}_{j=j_1}^{j_2}$ from
a $1-$Fr\'echet sample.

That is, the unknown parameter $\alpha$ appears only in the factor $1/\alpha^2$ of the covariance matrix and thus 
the GLS estimators $\what H_{\Sigma}$ and $\what C_{\Sigma}$ 
do not depend on $\alpha.$  Indeed, if one multiplies $\Sigma$ by a factor $\phi$, the resulting
estimates are not affected, since the formula \refeq{theta-GLS} involves the product of $\phi$ and its inverse.

 This invariance property shows that the GLS estimators can be computed {\it exactly}, 
 without using plug--in approximations for the unknown parameter $\alpha$ involved in
 the matrix $\Sigma$.   Table \ref{tab:psi} in the Appendix contains values of 
 $\psi(i)$ for $i=0,1\ldots,19$, obtained through Monte Carlo simulations.  This is sufficient to handle sample sizes of up to $
2^{20}=1, 048, 576$ observations.

\item Finally, $\Sigma_\alpha(j_1,j_2;N)$ is invertible, which follows from the fact that the joint 
distribution of the $Y_j$'s has a density with respect  to the Lebesgue measure. 
\end{enumerate}

\noi In view of the above remarks, we have that
\begin{corollary}
The minimum variance unbiased estimators for $H$ and $C$ in the regression model
\refeq{Yj-epsilonj},  linear in $Y_j$, are given by \refeq{theta-GLS}.  Moreover, the covariance matrix of $\what \theta_{\Sigma}$
 is 
$$
 \Sigma_{(\what H_{\Sigma},\what C_{\Sigma})} = (A^t \Sigma_\alpha^{-1}(j_1,j_2;N) A)^{-1} = \frac{1}{\alpha^2} (A^t \Sigma_1^{-1}(j_1,j_2;N) A)^{-1},
$$
where $\Sigma_1(j_1,j_2;N)$ is the covariance matrix of the $Y_j$ statistics based on $1-$Fr\'echet data.
\end{corollary}

\begin{figure}[ht!]
\begin{center}
\includegraphics[height=3in,width=4in]{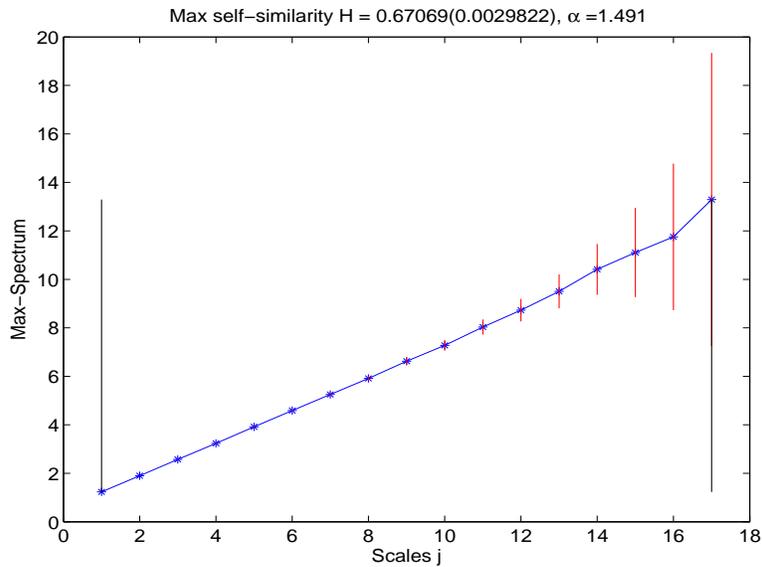}
{\caption{\label{fig:max_spectrum} \small Displayed is an example the max--spectrum of an i.i.d.\
$\alpha-$Fr\'echet sample of size $N = 2^{17} = 131, 072$ with $\alpha = 1.5$.  Observe that the max--spectrum is
 perfectly linear in $j$.  The vertical intervals around every $Y_j$ point indicate $95\%$ confidence intervals for
 the mean of $Y_j$ based on normal approximation.  Observe that these confidence intervals grow with the scale $j$. 
 GLS regression based on all scales $1\le j \le 17$ was used to obtain an estimate $\what \alpha = 1.491$.
 The estimated standard deviation of the slope $\what H = 0.67$ is indicated in parentheses: $\what \sigma_H = 0.00298$.
 This last estimate is based on the asymptotic variance of $\what H$ (see Proposition \ref{p:an-H-C}).}}
\end{center}
\end{figure}

In Figure \ref{fig:max_spectrum}, the max--spectrum of a sample from a Fr\'echet distribution with $N = 2^{17}$ observations 
is shown.  As expected, the max--spectrum is essentially linear in $j$ and the slope yields a very good
estimate of $1/\alpha$.  The asymptotic properties of estimators based on the max--spectrum of general heavy--tailed 
samples are established in Section \ref{s:cons-an}.  In practice, when the sample is not Fr\'echet, the max--spectrum is
linear in $j$ only on a range of the largest scales $j$.  The problem of choosing the ``best'' range of scales to estimate
$\alpha$ is very important in practice and is briefly addressed in Section \ref{s:cut-off}.

\section{Rates for moment--type functionals of heavy--tailed maxima}
\label{s:rate}

In this section, we establish some results for moment--type functionals obtained from maxima of
heavy--tailed data. They prove useful in establishing the consistency and asymptotic normality of the
max self--similar estimators under general conditions, but are also of independent interest since 
they yield {\em exact rates} of convergence in many cases.

 \medskip
 Let $X(1),X(2),\ldots,$ be i.i.d.\  random variables with c.d.f.\
 \beq\label{e:F-sigma}
  F(x) = \exp\{ -\sigma^\alpha(x) x^{-\alpha}\},\ \ x>0,
 \eeq
 where $\alpha>0$, and where the function $\sigma(x)>0$ is such that 
 $$
   \sigma(x)\longrightarrow \sigma_0>0,\ \ \mbox{ as }x\to\infty.
 $$ 
Here, we let the function $\sigma(x)$ take values in the extended half--line $(0,\infty]$, that is,
$\sigma(x)$ can take the value $\infty$, in which case $F(x)$ becomes $e^{-\infty} = 0$ (see the
Examples below).  Such a representation always exists if the c.d.f.\  $F$ belongs to the normal
domain of attraction of an $\alpha-$Fr\'echet distribution, that is, if
\beq\label{e:Mn-Z}
  M_n := \frac{1}{n^{1/\alpha}} \bigvee_{1\le i \le n} X(i) \stackrel{d}{\longrightarrow} Z,
\eeq
where $G(x) := \P\{Z \le x\} = \exp\{ - \sigma_0^\alpha x^{-\alpha}\},\ x>0$, for some $\sigma_0>0$.
For simplicity, we suppose that the $X(i)$'s are positive, almost surely, that is $F(0)=0$.
The case when the $X(i)$'s can take negative values is addressed in Section \ref{s:cons-an} below.
  
Our goal here is to establish bounds on the rate of convergence of $\E f(M_n)$ to $\E f(Z)$, as
$n\to\infty$, for an absolutely continuous function $f:(0,\infty)\to\bbR$.  We do so under general
conditions on the asymptotic tail behavior of the c.d.f.\ $F(x)$.  

\medskip
\noi In what follows, the next two conditions on the c.d.f.\  $F(x)$ are needed:

\begin{condition} For some $\beta>0$ and $C_1>0$,
\label{cond:C1} 
\beq\label{e:C1}
 |\sigma^\alpha(x) - \sigma_0^\alpha| \le C_1 x^{-\beta},\ \mbox{ for all sufficiently large }x>0.
\eeq
\end{condition}

\noi and

\begin{condition}
\label{cond:C2} We have $F(0) = 0$ and for some $C_2>0$,
 \beq\label{e:C2}
  \sigma^\alpha(x)  \ge C_2 \min\{ 1, x^\gamma\} ,\ x>0, \ \mbox{ for some }\ \gamma \in (0,\alpha).
\eeq
\end{condition}

 In the examples below, we show that the Conditions \ref{cond:C1} and {\ref{cond:C2} hold in many cases of
 practical interest.  The second condition concerns the behavior of $F(x)$ for small $x$, and ensures that
 $\E(X^p 1_{\{X\le 1\}}) <\infty$, for any $p\in\bbR$.  This condition always holds, for example, if the $X(i)$'s
 are bounded away from zero, almost surely.  The case of arbitrary $X(i)$'s which can possibly take negative
 values is addressed in Section \ref{s:cons-an}.

\medskip
The following result provides an upper bound on $|\E f(M_n) - \E f(Z)|$ under the above conditions for
general class of absolutely continuous functions $f$.  Namely, we shall suppose that
$f(x) = f(x_0) + \int_{x_0}^x f'(u) du,\ x>0$, for some (any) $x_0\in(0,\infty)$, 
with $f'$ being a locally integrable function.

\begin{theorem}\label{t:rate}
 Let $f(x), x>0$ be an absolutely continuous function on all compact intervals $[a,b] \subset (0,\infty)$.  
 Let also $F_n(x):= \P\{M_n \le x\}$ and $G(x)= \P\{Z\le x\},\ x\in\bbR$, be the c.d.f.'s of the random variables
 $M_n$ and $Z$ in \refeq{Mn-Z}.  Suppose that Conditions \ref{cond:C1} and \ref{cond:C2} hold.
 
{\bf (a)} If for some $m \in \bbR$ and $\delta>0$,
 \beq\label{e:p:rate-a}
  x^m|f(x)| + \mathop{{\rm esssup}}_{0 < y \le x} y^m|f'(y)| \to 0,\ x\downarrow 0,\ \ \mbox{ and }
  \ \ x^{-\alpha}|f(x)| + x^{1+\delta} \mathop{{\rm esssup}}_{y\ge x} y^{-\alpha}|f'(y)| \to 0,\ x\to\infty,
 \eeq
then $\E |f(Z)|$ and $\E |f(M_n)|,\ n\in \bbN$ are finite.  Moreover,
\beq\label{e:p:rate-a-1}
 \E f(M_n) - \E f(Z) = \int_0^\infty (G(x) - F_n(x)) f'(x) dx. 
\eeq
Here ${\rm esssup}$ denotes the essential supremum of a measurable function $g$, that is, 
$$
  {\rm esssup}_{y\in A} g(y) := \inf_{A_0 \subset A,\ |A\setminus A_0| =0} \sup_{y\in A_0} g(y),
$$
for any Borel set $A$, where $|A|$ denotes the Lebesgue measure of the set $A$.

{\bf (b)} If in addition to \refeq{p:rate-a}, $\int_1^\infty x^{-(\alpha+\beta)} |f'(x)| dx < \infty$, then
 for any $\epsilon(n) \to 0,$ such that $n^{1/\alpha}\epsilon(n)\to \infty,$ as $n\to\infty$, we have
\begin{eqnarray}\label{e:p:rate-b}
 |\E f(M_n) - \E f(Z)| &\le&  C_1 n^{-\beta/\alpha} {\Big(} \int_0^\infty  x^{-(\alpha+\beta)}
 |f'(x)| e^{-c x^{-\alpha}} dx {\Big)} \nonumber \\
 & & \ \ \ \ \ \ + 2 \int_0^{\epsilon(n)} e^{-C_2 x^{-(\alpha-\gamma)}} | f'(x)| dx,
\end{eqnarray}
for all sufficiently large $n$, where $c \in (0,\sigma_0^\alpha)$ can be chosen arbitrarily close
to $\sigma_0^\alpha$.  Moreover, 
\beq\label{e:p:rate-b-order}
 |\E f(M_n) - \E f(Z)| \le  C_f n^{-\beta/\alpha},
\eeq
for all sufficiently large $n$ with some $C_f>0$.
\end{theorem}
\begin{proof}
{\it We first prove part {\bf (a)}.}  Let $f(x) = f(x_0) + \int_{x_0}^x f'(u)du,\ x>0,$ with $x_0\in (0,\infty)$,
where $f'(x),\ x\in(0,\infty)$ is locally integrable, and where $\int_a^b = - \int_{b}^a$. 
Let now $[a,b] \subset (0,\infty),\ x_0\in(a,b)$ be an arbitrary interval
and observe that $\int_{a}^b f(x) d F_n(x)$ equals
\begin{eqnarray}
 & &\int_a^{x_0} f(x) d F_n(x) + \int_{x_0}^b f(x) d (F_n(x) -1) 
 = F_n(x_0) f(x_0) - F_n(a) f(a) - \int_a^{x_0} F_n(x) f'(x) dx \nonumber\\
 & & \ \ \ \ \  \ \ \ + (F_n(b)-1)f(b) - (F_n(x_0)-1)f(x_0) -
 \int_{x_0}^b (F_n(x) - 1)f'(x) dx\label{e:p:rate-1} \\ 
 & & = (F_n(b)-1)f(b) - F_n(a)f(a) + f(x_0)  \nonumber\\
 & & \ \ \ \ \ \ \ \ \  - \int_a^{x_0} F_n(x) f'(x) dx + \int_{x_0}^b (1-F_n(x))f'(x)dx.  \label{e:p:rate-2}
\end{eqnarray}
The equality in Relation \refeq{p:rate-1} follows from Lemma \ref{l:int-by-parts}.

In view of Relation \refeq{p:rate-2}, the monotone convergence theorem implies that
$\E |f(M_n)| = \int_0^\infty |f(x)| d F_n(x)$ is finite if
\beq\label{e:p:rate-1.25} 
 |(F_n(b)-1)f(b)| + |F_n(a)f(a)| \longrightarrow 0,\ \mbox{ as }a\downarrow 0\ \ \mbox{  and  }\ \ b\to\infty,
\eeq
{\it and} if
\beq\label{e:p:rate-1.5}
 \int_0^{x_0} F_n(x)|f'(x)|dx + \int_{x_0}^\infty (1-F_n(x)) |f'(x)| dx  < \infty.
\eeq
Observe that by \refeq{F-sigma},
$$
 F_n(x) = F(n^{1/\alpha}x)^n = \exp\{-\sigma^\alpha(n^{1/\alpha}x) x^{-\alpha}\},\ \ x>0.
$$
Hence, in view of \refeq{C1} we have
\beq\label{e:p:rate-3}
 1-F_n(x) \sim \sigma_0^\alpha x^{-\alpha},\ \mbox{ as }x\to\infty,
\eeq 
since $1- e^{-u} \sim u,$ as $u\downarrow 0$.
Thus, the second convergence in \refeq{p:rate-a}, implies $|(F_n(b)-1)f(b)| \to 0,\ b\to\infty$.
On the other hand, by \refeq{C2}, for $n\ge 1,\ n\in\bbN$,
\beq\label{e:p:rate-4}
 \sigma^\alpha(n^{1/\alpha}x) \ge C_2 n^{\gamma/\alpha} x^{\gamma} \ge C_2 x^{\gamma},\ \mbox{ for all }
 x\in (0,n^{-1/\alpha}),
\eeq
and hence 
\beq\label{e:p:rate-3.1-new}
  F_n(x) = \exp\{-\sigma^\alpha(n^{1/\alpha}x)x^{-\alpha}\} \le \exp\{-C_2 x^{-(\alpha-\gamma)}\},\ \ 
   \mbox{ for all }\ x\in(0,n^{-1/\alpha}).
\eeq
Thus, since $u^pe^{-u} \to 0,$ as $u\to\infty$, for any $p\in\bbR$, the first convergence
in \refeq{p:rate-a} implies that $F_n(a) f(a) \to 0$, as $a\to\infty$.  We have thus shown that
\refeq{p:rate-1.25}  holds.  One can similarly show that the integrals in \refeq{p:rate-1.5} are finite by the
 using the conditions in \refeq{p:rate-a} on $f'$ and Relations \refeq{p:rate-3} and \refeq{p:rate-4}.  Indeed, 
for almost all $x >0$, we have
\beq\label{e:p:rate-4.5}
  F_n(x) |f'(x)| \le ({\rm sup}_{0<y\le x} F_n(y)y^{-m}) ({\rm esssup}_{0<y\le x} y^{m}|f'(y)|) =
   {\cal O}(x^{-|m|} \exp\{-C_2 x^{-(\alpha-\gamma)}\}) \to 0,
\eeq
as $x\downarrow 0$ and, for almost all $x>0$,
\beq\label{e:p:rate-4.6}
(1-F_n(x)) |f'(x)| \le  (\sup_{y\ge x} (1-F_n(y))y^{-\alpha}) ({\rm esssup}_{y\ge x} y^{\alpha} |f'(y)|)
= {\cal O}(x^{-(1+\delta)}),
\eeq
as $x\to\infty.$
We have thus shown that $\int_0^\infty|f(x)|dF_n(x) <\infty$ for all $n\in\bbN$.
One can similarly show that $\int_0^\infty|f(x)|dG(x) <\infty$, by replacing $F_n(x)$ with $G(x)$, above, and using
the fact that $G(x) = \exp\{-\sigma_0^\alpha x^{-\alpha}\},\ x>0$ satisfies trivially Conditions \ref{cond:C1} 
and \ref{cond:C2}.

Observe that \refeq{p:rate-a-1} follows from the relations
$$\int_0^{\infty} f(x) d F_n(x) = f(x_0) - \int_0^{x_0} F_n(x) f'(x)dx + \int_{x_0}^\infty (1-F_n(x))f'(x) dx$$ 
and 
$$\int_0^{\infty} f(x) d G(x) = f(x_0) - \int_0^{x_0} G(x) f'(x)dx + \int_{x_0}^\infty (1-G(x))f'(x) dx.$$

\medskip
{\it We now turn to proving part {\bf (b)}.}  Let $\epsilon(n)\downarrow 0$ be such that $n^{1/\alpha}\epsilon(n)\to
\infty$, as $n\to\infty$.  By \refeq{p:rate-a-1}, using the triangle inequality, we get
\begin{eqnarray*}
 |\E f(M_n) - \E f(Z)| & \le& \int_0^{\epsilon(n)} G(x) |f'(x)| dx + \int_0^{\epsilon(n)} F_n(x) |f'(x)| dx \\
 & & \ \ \ \ \ \ \ +  \int_{\epsilon(n)}^\infty |F_n(x) - G(x)| |f'(x)| dx =:  I_1 + I_2 + I_3.
\end{eqnarray*}
We first consider the integral $I_3.$  Since $n^{1/\alpha}\epsilon(n) \to \infty,\ n\to\infty$, in view of 
\refeq{C1}, for all sufficiently large $n$, we have 
\begin{eqnarray}\label{e:p:rate-5}
 |F_n(x) - G(x)| & = & |\sigma^\alpha(n^{1/\alpha}x) - \sigma_0^\alpha| x^{-\alpha} e^{-\theta_n(x) x^{-\alpha}}\nonumber\\
                 &\le& C_1 n^{-\beta/\alpha} x^{-(\alpha+\beta)} e^{-c x^{-\alpha}},
\end{eqnarray}
for all $x \in (\epsilon(n),\infty)$, where $c$ is an arbitrary constant in $(0,\sigma_0^\alpha)$, and where $\theta_n(x)$
is between $\sigma^\alpha(n^{1/\alpha}x)$ and $\sigma_0^\alpha$.  Indeed, the first relation in \refeq{p:rate-5} follows
by the mean value theorem applied to the function $g(u) = \exp\{ - u x^{-\alpha}\},\ u>0$.  The inequality in \refeq{p:rate-5},
follows from \refeq{C1} since $n^{1/\alpha}\epsilon(n) \to\infty$ implies
$\sup_{x\ge \epsilon(n)} \sigma^\alpha(n^{1/\alpha}x) \ge c,\ c\in(0,\sigma_0^\alpha),$
for all sufficiently large $n$.

Therefore \refeq{p:rate-5} implies
$$
 I_3 \le C_1 n^{-\beta/\alpha} \int_{\epsilon(n)}^\infty x^{-(\alpha+\beta)} e^{-c x^{-\alpha}} |f'(x)| dx
 \le C_1 n^{-\beta/\alpha} \int_0^\infty x^{-(\alpha+\beta)} |f'(x)|  e^{-c x^{-\alpha}} dx,
$$
for all sufficiently large $n$.   The last integral is finite.  Indeed, by assumption
$\int_1^{\infty} x^{-(\alpha+\beta)}  |f'(x)| dx < \infty$.  The integral $\int_0^{1}
x^{-(\alpha+\beta)} |f'(x)| e^{-cx^{-\alpha}} dx $ is finite 
since in view of \refeq{p:rate-a},
\beq\label{e:p:rate-5.1-new}
 ({\rm esssup}_{0\le y\le x} y^{m} |f'(y)|) x^{-(\alpha+\beta+|m|)}  
 e^{-c x^{-\alpha}} = {\cal O}{\Big(} x^{-(\alpha+\beta+|m|)} e^{-cx^{-\alpha}} {\Big)} 
 = {\cal O}(x^p),\ x\downarrow 0,
\eeq
for any $p>0$.

\smallskip
{\it We now consider the integral $I_2$.} Observe that $\epsilon(n) > n^{-1/\alpha},$ eventually, and hence
\beq\label{e:p:rate-6}
I_2 \le \int_0^{n^{-1/\alpha}} \exp\{ - C_2 x^{-(\alpha-\gamma)}\} |f'(x)| dx + \int_{n^{-1/\alpha}}^{\epsilon(n)}
F_n(x) |f'(x)| dx,
\eeq
by \refeq{p:rate-3.1-new}. Relation \refeq{C2} implies that $\sigma^\alpha(n^{1/\alpha}x) \ge C_2,$ for all
$x\in(n^{-1/\alpha},\epsilon(n))$, and hence
$F_n(x) \le \exp\{ - C_2 x^{-\alpha}\} \le  \exp\{ - C_2 x^{-(\alpha-\gamma)}\}$, $x\in(n^{-1/\alpha},\epsilon(n))$.
Therefore, the second integral in \refeq{p:rate-6} can be bounded above by $\int_{n^{-1/\alpha}}^{\epsilon(n)}
 \exp\{-C_2x^{-(\alpha-\gamma)}\} |f'(x)| dx$ and hence
$$
I_2 \le \int_0^{\epsilon(n)} \exp\{-C_2x^{-(\alpha-\gamma)}\} |f'(x)| dx.
$$
One can similarly bound $I_1$.  Indeed, Relation \refeq{C2} implies that $\sigma_0^\alpha \ge C_2$, 
since $\sigma^\alpha (x)\sim \sigma_0^\alpha,\ x\to\infty$.  For all
$0<x<\epsilon(n)<1$ and $\gamma \in(0,\alpha)$, we have $x^{-\alpha} \ge x^{-(\alpha-\gamma)},$ and hence
we obtain
$$
  I_1 = \int_{0}^{\epsilon(n)} \exp\{ -\sigma_0^\alpha x^{-\alpha}\} |f'(x)| dx 
      \le \int_0^{\epsilon(n)} \exp\{-C_2x^{-(\alpha-\gamma)}\} |f'(x)| dx.
$$
The last three bounds for $I_1$, $I_2$ and $I_3$ imply \refeq{p:rate-b}.  

{\it Now, to prove \refeq{p:rate-b-order}}, observe that, as in \refeq{p:rate-5.1-new}, 
since $\alpha-\gamma>0$, for almost all $x>0$, we have
\beq\label{e:p:rate-7}
\exp\{-C_2x^{-(\alpha-\gamma)}\} |f'(x)| \le {\cal O}{\Big(} x^{-|m|} e^{-C_2x^{-(\alpha-\gamma)}} {\Big)}
 = {\cal O}(x^p),\ x\downarrow 0,
\eeq
for any $p>0$.  Thus, the second integral in \refeq{p:rate-b} is of order ${\cal O}(\epsilon(n)^p)$, for any
$p>0$ and by setting $\epsilon(n):= n^{-\delta},$ for some $\delta \in (0,1/\alpha)$, we obtain that 
\refeq{p:rate-b-order} holds. This completes the proof of the theorem.
$\Box$
\end{proof}

\medskip
In the following examples we show that most heavy--tailed distributions of practical interest
satisfy the conditions of Theorem \ref{t:rate}.

 \medskip
 \noi{\bf Examples:}
\begin{itemize}
 \item {\it (Pareto laws)} Let $F(x) = 1- (x/\sigma_0)^{-\alpha},\ x\ge \sigma_0,$ 
 and $F(x)= 0,$ $x<\sigma_0$, for some $\sigma_0>0$ and $\alpha>0$.  Then, Relation
 \refeq{F-sigma} holds with
 $$
  \sigma^\alpha(x) = \infty 1_{(0,\sigma_0]}(x) 
                    - x^{\alpha}\ln(1-(x/\sigma_0)^{-\alpha})1_{(\sigma_0,\infty)}(x),
 $$
 that is, the function $\sigma(x)$ equals $\infty$ for all $x\in(0,\sigma_0]$ to 
 account for the fact that $F(x) = 0,\ x\in(0,\sigma_0]$.

{\it Observe that $\sigma^\alpha(x)$ satisfies Condition \ref{cond:C1} with 
 $\beta=\alpha.$}  Indeed, since $\ln(1-u) = -u + u^2/2 + {\cal O}(u^3),\ u\to 0$, by
 setting $u:= (x/\sigma_0)^{-\alpha}$, 
 we obtain
 \beq\label{e:pareto-C1}
  |\sigma^\alpha(x) -\sigma_0^\alpha| =
  {\Big|} { \ln(1-(x/\sigma_0)^{-\alpha}) \over x^{-\alpha}} + \sigma_0^\alpha {\Big|}
  = \sigma_0^\alpha {\Big|} {\ln(1-u)\over u} + 1{\Big|} \le \sigma_0^\alpha u
  = \sigma_0^{2\alpha} x^{-\alpha}, 
 \eeq
 for all sufficiently large $x$.

 One has, moreover, that
 \beq\label{e:Pareto-exact-rate}
   \sigma^\alpha(x) -\sigma_0^\alpha \sim \frac{\sigma_0^{2\alpha}}{2} x^{-\alpha},\ \ 
    \mbox{ as }x\to\infty.
 \eeq 
 (see Proposition \ref{p:rate-exact}, below).

{\it Condition \ref{cond:C2} also holds.}  Indeed, $\sigma(x)=\infty \ge x^\gamma$, for
 all $x \in (0,\sigma_0]$ and $\gamma \in (0,\alpha)$.  To prove \refeq{C2}, it remains
 to show that $\sigma^\alpha(x) \ge C_2  >0$, for all $x>0$.  As shown in \refeq{pareto-C1}
 above $\sigma^\alpha(x) \to \sigma_0^\alpha,\ x\to\infty$, where $\sigma_0>0$.  On the 
 other hand $\sigma^\alpha(x)$ is a positive, continuous function over all compact 
 intervals of $(\sigma_0,\infty)$ and $\sigma(x) \to \infty$, as $x\to\sigma_0$.  This 
 shows that $\sigma^\alpha(x)$ is bounded below by a positive constant.

 \item {\it (Products of Fr\'echet laws)} Let $F(x) = G_{\alpha_0}(x/\sigma_0)
 G_{\alpha_1}(x/\sigma_1),$ where $\sigma_0,\ \sigma_1 >0$ and $0<\alpha_0<\alpha_1$, and where 
 $G_\alpha(x) = \exp\{-x^{-\alpha}\},\ x>0$ denotes the c.d.f.\ of a standard $\alpha-$Fr\'echet
 variable.  Observe that the function $F(x)$ is the c.d.f.\  of $\max\{\sigma_0 Z_0, \sigma_1 Z_1\},$
 where $Z_0$ and $Z_1$ are independent standard $\alpha_0-$ and $\alpha_1-$Fr\'echet random
 variables, respectively.  Therefore, \refeq{F-sigma} holds with $\alpha=\alpha_0$ and
 \beq\label{e:product-Frechet-exact}
   \sigma^\alpha(x) = \sigma_0^\alpha + \sigma_1^\alpha x^{-(\alpha_1-\alpha_0)},\ \ x>0.
 \eeq
 Conditions \ref{cond:C1} and \ref{cond:C2} are readily satisfied where $\beta = \alpha_1-\alpha_0>0$.

 \item {\it (Mixtures of Pareto laws)} Let 
 $$
  F(x) = p (1 - (x/\sigma_0)^{-\alpha_0})1_{\{x\ge \sigma_0\}} + 
   (1-p) (1-(x/\sigma_1)^{-\alpha_1})1_{\{x\ge \sigma_1\}},\ \ 0<\alpha_0<\alpha_1,
 $$
 where $p\in(0,1)$ and $\sigma_0,\ \sigma_1>0$.
 
 Then, \refeq{F-sigma} holds with $\alpha\equiv \alpha_0$, and
 $
 \sigma^\alpha(x) = \infty 1_{(0,\sigma_*]}(x) - x^{\alpha} \ln(F(x)) 1_{(\sigma_*,\infty)}(x),
 $
 where $\sigma_* := \min\{\sigma_0,\sigma_1\} >0$.

 As in the case of Pareto laws, one can show that Condition \ref{cond:C1} holds with $\beta = 
 \min\{\alpha_0, \alpha_1-\alpha_0\}$ and, $\sigma_0$ replaced by $p\sigma_0$.  In fact,
 \beq\label{e:mix-Pareto-exact}
  \sigma^\alpha(x) - p\sigma_0^\alpha \sim C_0 x^{-\beta},\ \ \mbox{ as } x\to\infty,
 \eeq
 where 
 $$
 C_0 = \left\{\begin{array}{ll}
              \sigma_1^{\alpha_1}(1-p) &,\ \mbox{ if }\alpha_1-\alpha_0 < \alpha_0\\ 
               \sigma_1^{\alpha_1}(1-p) + p^2\sigma_0^{2\alpha_0}/2 &,\ \mbox{ if }\alpha_1-\alpha_0 = \alpha_0\\
                 p^2\sigma_0^{2\alpha_0}/2 &,\ \mbox{ if }\alpha_1-\alpha_0 > \alpha_0    
             \end{array}
 \right.
 $$
 One can also show that Condition \ref{cond:C2} holds as in the case of Pareto laws. 
 
 \item Absolute values of $\alpha-$stable ($0<\alpha<2$) and $t-$distributed random variables $X_i$'s, for example, 
 also satisfy Condition \ref{cond:C1}.  They {\it do not} satisfy  Condition \ref{cond:C2}, however, since 
 $\E (|X_1|^{-1} 1_{\{|X|\le 1\}})$ is infinite.  In Proposition \ref{p:coupling} below, we address the general case 
 where Condition \ref{cond:C2} fails and in fact the case where the $X_i$'s can take negative values.
\end{itemize}

\medskip
\noi The following result shows that the rate $n^{-\beta/\alpha}$ in \refeq{p:rate-b-order} is optimal, if
so is the inequality in \refeq{C1}.  

\begin{proposition}\label{p:rate-exact}   Assume that $F$ is as in \refeq{F-sigma} and satisfies Conditions \ref{cond:C1} and 
 \ref{cond:C2} above, and let $f$ be as in Theorem \ref{t:rate} {\bf (b)}.  Suppose, in addition, that
 $\sigma^\alpha(x) - \sigma_0^\alpha \sim C_1 x^{-\beta},$ as $x\to\infty$, for some $C_1\not=0$.  Then
 \beq\label{e:p:rate-exact}
  n^{-\beta/\alpha} (\E f(M_n) - \E f(Z)) \longrightarrow
   C_1 \int_0^{\infty} x^{-(\alpha+\beta)} f'(x) e^{-\sigma_0^\alpha x^{-\alpha}} dx,
  \ \mbox{ as } n\to\infty. 
 \eeq
 \end{proposition}
\begin{proof}
  Let as in Theorem \ref{t:rate}, $\epsilon(n)\to 0$ be such that
  $n^{1/\alpha}\epsilon(n)\to\infty,$ as $n\to\infty$.  The triangle inequality
  applied to Relation \refeq{p:rate-a-1} implies
  \beq\label{e:p:rate-exact-1}
  {\Big|} \E f(M_n) - \E f(Z)  - \int_{\epsilon(n)}^\infty (G(x)-F_n(x))f'(x) dx {\Big|}
  \le \int_{0}^{\epsilon(n)} G(x) |f'(x)| dx + \int_{0}^{\epsilon(n)} F_n(x) |f'(x)| dx.
  \eeq
  As in the proof of Theorem \ref{t:rate} one can show that the integrals in the right--hand side of the 
  last expression are of order $o(n^{-\beta/\alpha}),$ as $n\to\infty$, if $\epsilon(n) := n^{-\delta},$ 
  $\delta\in (0,1/\alpha)$ (see \refeq{p:rate-7}).

  To establish \refeq{p:rate-exact} we will now examine the order of the integral in
  the left--hand side of \refeq{p:rate-exact-1}.
  Observe that
 \beq\label{e:p:rate-exact-2}
  { \sigma^\alpha(n^{1/\alpha}x) - \sigma_0^\alpha \over n^{-\beta/\alpha}} \longrightarrow
  C_1 x^{-\beta},
 \eeq 
 as $n\to\infty$, for all $x>0$.
 Hence (as in Theorem \ref{t:rate}), in view of \refeq{F-sigma} and \refeq{p:rate-exact-2},
 the mean value theorem implies
 $$
  {n^{\beta/\alpha}}{(G(x)-F_n(x))f'(x)} \longrightarrow C_1 x^{-(\alpha+\beta)}f'(x) e^{-\sigma_0^\alpha x^{-\alpha}}, 
 $$ 
 as $n\to\infty$, for any $x\in(\epsilon(n),\infty)$ and hence for any $x>0$ 
 ($\epsilon(n)\to 0$, $n\to\infty$).  
 As in the proof of Theorem \ref{t:rate}, one can show that the left--hand side of the
 last expression is bounded above in absolute value by an integrable function.  Therefore, the dominated 
 convergence theorem implies that $n^{\beta/\alpha}\int_{0}^{\infty}(G(x)-F_n(x))f'(x) dx$ converges to
 the integral in \refeq{p:rate-exact}, as $n\to\infty$. 
 $\Box$
\end{proof}

\medskip
\noi The next result, which follows directly from Theorem \ref{t:rate} is used in 
Section \ref{s:cons-an}.

 \begin{corollary}\label{c:logs} Assume that $F$ is as in \refeq{F-sigma} and satisfies Conditions \ref{cond:C1} and 
 \ref{cond:C2} above.  Then $\E |\ln(M_n)|^p <\infty$ for all $n\in\bbN$ and $p>0$.  Moreover,
 for any $p>0$ and $k\in\bbN$, we have
 $$
    {\Big|} \E |\ln(M_n)|^p - \E |\ln(Z)|^p {\Big|} = {\cal O}(n^{-\beta/\alpha})\  \ \mbox{ and }\ \
   {\Big|}\E \ln(M_n)^k - \E \ln(Z)^k {\Big|} = {\cal O}(n^{-\beta/\alpha}),
 $$
 as $n\to\infty$, where $M_n$ and $Z$ are as in Theorem \ref{t:rate}.
 \end{corollary}

\medskip
 In Section \ref{s:cons-an}, one encounters covariance functionals of maxima over blocks of 
 heavy--tailed variables, that is, bivariate moment--type functionals arise.  The following
 result establishes rates of convergence for such functionals in the special case of logarithms.

 \begin{corollary}\label{c:cov-logs}
 Suppose that $F$ is as in \refeq{F-sigma} and satisfies Conditions \ref{cond:C1} and 
 \ref{cond:C2}.  Let $X(1),\ldots,X(n)$ and $Y(1),\ldots, Y(m),\ n,m\in\bbN$
 be i.i.d.\  random variables with c.d.f.\  $F(x).$  Consider the normalized maxima
 $$
 M_n^X := \frac{1}{n^{1/\alpha}}\bigvee_{1\le i \le n} X(i) \ \ \mbox{ and }\ \ 
 M_m^Y := \frac{1}{m^{1/\alpha}}\bigvee_{1\le i \le m} Y(i),\ n, m\in\bbN.
 $$
 Then, for any $a>0$, as $n,\ m\to\infty$, we have that
 \beq\label{e:c:cov}
  \E \ln(M_n^{X}) \ln(M_n^X \vee aM_m^Y) - \E \ln(Z_X)\ln(Z_X \vee a Z_Y) = {\cal O}(n^{-\beta/\alpha} 
 +m^{-\beta/\alpha}),
 \eeq
 where $Z_X$ and $Z_Y$ are independent $\alpha-$Fr\'echet random variables with scale coefficients
 $\sigma_0$. 
 \end{corollary}

 Corollary \ref{c:cov-logs} was stated  in generality which allows us to have different number of
 $X(i)$'s and  $Y(i)$'s ($n$ and $m$, respectively) in the maxima $M_n^X$ and $M_m^Y$.  This 
 flexibility is needed for the proof of Proposition \ref{p:mom-cov} below.

\medskip
\noi{\sc Proof of Corollary \ref{c:logs}:}
 Let $f(x) = |\ln(x)|^p,\ p>0,\ x>0$.  Observe that $f(x) = \int_1^x f'(u) du,$
 where $f'(x) = p|\ln(x)|^{p-1}/x$ for $x\ge 1$ and $f'(x) = -p|\ln(x)|^{p-1}/x,$ for $0<x\le 1$.
 One can verify that the conditions in \refeq{p:rate-a} are fulfilled and therefore, Theorem
 \ref{t:rate} implies the result.  The argument in the case when $f(x) = (\ln(x))^k,\ k\in\bbN$ is
 similar. $\Box$

\medskip
\noi{\sc Proof of Corollary \ref{c:cov-logs}:}
 By Corollary \ref{c:logs}, the expected values in \refeq{c:cov} exist
 since $\E |\ln(M_n^X)|^p <\infty,\ \forall p>0$ and since $a\vee b \le a + b$ for any $a,\ b\ge 0$.
 Observe that by independence and Fubini's theorem,
 $$
  \E \ln(M_n^{X}) \ln(M_n^X \vee aM_m^Y) = \int_0^\infty {\Big(} \int_0^\infty f(x,y) d F_n(x){\Big)} d F_m(y), 
 $$
 and 
 $$
  \E \ln(Z_X)\ln(Z_X\vee a Z_Y) = \int_0^\infty  {\Big(} \int_0^\infty f(x,y) d G(x){\Big)} d G(y), 
 $$
 where $f(x,y) = \ln(x)\ln(x\vee ay),\ x,y,a>0$, $F_n(x):= F(n^{1/\alpha}x)^n$ is the
 c.d.f.\  of $M_n^X$ (and $M_n^Y$), and where $G(x) = \exp\{-\sigma_0^\alpha x^{-\alpha}\},\ x>0$.
 Now, by adding and subtracting the term $\int_0^\infty (\int_0^\infty f(x,y) d G(x))d F_m(y)$, applying
 Fubini's theorem and then the triangle inequality, we obtain that the left--hand side of \refeq{c:cov} is
 bounded above in absolute value by
 \begin{eqnarray*}
 & &  \int_0^\infty {\Big|}\int_0^\infty f(x,y) dF_n(x) - \int_0^\infty f(x,y) d G(x){\Big|} dF_m(y) \\
 & & \ \ \ \ \ \ \ \ + \int_0^\infty {\Big|}\int_0^\infty f(x,y) dF_m(y) - \int_0^\infty f(x,y) d G(y) {\Big|}
  dG(x) =: I_1 +I_2.
 \end{eqnarray*}
 Focus next on the term $I_1$.  Let $g(y):= \int_0^\infty f(x,y) (dG(x) - dF_n(x)), y>0$.
 Observe that for each $y>0,\ y\not=x/a$, $f(x,y)$ is differentiable in $x$ since
 $$ f(x,y) =\left\{\begin{array}{ll}
                     \ln(x)\ln(ay) &,\ 0<x<ay \\
                     \ln(x)^2 &,\ ay \le x
                   \end{array}
             \right.
 $$
In fact,
$$
  |f'_x(x,y)| \le 2|\ln(x)|/x + |\ln(ay)|/x,\ \ x>0,\ y>0.
$$
Thus, Theorem \ref{t:rate} {\bf (b)}, applied to the inner integral $g(y)$ in $I_1$ implies
 \beq\label{e:c:cov-2}
  |g(y)| \le n^{-\beta/\alpha}(C' + C''|\ln(y)|),
 \eeq
for all sufficiently large $n$, where the constants $C'>0$ and $C''>0$ do not depend on $y$
(This follows from Relation \refeq{p:rate-b} by taking $\epsilon(n):= n^{-\delta},\ \delta\in(0,1/\alpha)$
 and observing that the second integral therein is negligible with respect to the term 
 $(1+|\ln(y)|)n^{-\beta/\alpha}$.)  

 Note now that the function $|\ln(y)|$ satisfies the assumptions of Theorem \ref{t:rate} {\bf (b)} 
 and hence $\int_0^\infty |\ln(y)| dF_m(y) \to \int_0^{\infty} |\ln(y)|dG(y)$,
 as $m\to\infty$.  Therefore, the inequality \refeq{c:cov-2} implies that
 $I_1 = {\cal O} (n^{-\beta/\alpha}),$ as $n\to\infty$.  One can similarly show that 
 $I_2 = {\cal O}(m^{-\beta/\alpha}),\ m\to\infty$.  
 $\Box$

\section{Asymptotic properties of the max self--similarity estimators}
 \label{s:cons-an}

 We establish here the consistency and asymptotic normality of the
 estimators defined in \refeq{H-hat}, above.  In fact, we prove joint asymptotic
 normality of the max self--similarity estimators of the tail exponent $\alpha$
 and the scale coefficient $\sigma_0$.  These results rely on the behavior of
 moment--type functionals of heavy--tailed maxima established in Section \ref{s:rate}.

 The general case
 where the $X(i)$'s may be $0$ or even take negative values is addressed at the end of this section.

\medskip
Let the $Y_j$'s be defined as in \refeq{Yj}, where now $N$ denotes the sample size of
available $X(i)$'s, $1\le j \le [\log_2 N]$ and where $N_j:= [N/2^j]$.
As noted above, the larger the scales $j$, the more precise the asymptotic relation \refeq{Yj-sim}.
Therefore, to obtain consistent estimates
for the parameter $H=1/\alpha$ one should focus on a range of scales which grows as
the sample size increases.  We therefore {\it fix} a range $j_1\le j \le j_2,\ j_1, j_2\in \bbN$
and focus on the vectors 
$$
 Y_r:= \{ Y_{j+r}\}_{j=j_1}^{j_2},
$$ 
with $r \in \bbN,\ j_2+r\le [\log_2 N]$ where the parameter $r=r(N)$ grows
with the sample size.

 The following result shows that the mean and the covariance matrix of the vector $\bY_r$
 are asymptotically equivalent to the mean and and the covariance matrix 
 in the case where the $X(i)$'s are $\alpha-$Fr\'echet (see Proposition \ref{p:cov}).

 \begin{proposition}\label{p:mom-cov} Suppose that the c.d.f.\  $F$ has the representation \refeq{F-sigma}
 and satisfies Conditions \ref{cond:C1} and \ref{cond:C2}, above.   

 Then, 
 \beq\label{e:moments}
   {\Big|} \E Y_{j+r} - \mu_r(j) {\Big|} = {\cal O}{\Big(} 1/2^{r\beta/\alpha}{\Big)},
  \ \ \mbox{ as } r\to\infty,
 \eeq
 and for any fixed $j_1\le i\le j\le j_2,\ i,j\in\bbN$, we have
 \beq\label{e:cov}
  {\Big|} N_{j_2+r} {\rm Cov}(Y_{i+r},Y_{j+r}) - \alpha^{-2} \Sigma_1(i,j) {\Big|}
   = {\cal O}{\Big(} 1/2^{r\beta/\alpha}{\Big)} + {\cal O}{\Big(}2^r/N{\Big)},\ \ \mbox{ as } r\to\infty.
 \eeq
  Here 
 \beq\label{e:mu_r-Sigma} 
   \mu_r(j):= (j+r)/\alpha + C(\sigma_0,\alpha)\ \ \mbox{ and }\ \ 
   \Sigma_1(i,j) = 2^{j-j_2} \psi(|i-j|),
 \eeq
 where the function $\psi$ is defined in \refeq{psi} and where $C(\sigma_0,\alpha)$ is as in \refeq{the-C}.
 \end{proposition}
 \begin{proof}  Observe that by \refeq{Yj}, we have
 $
  \E Y_{j+r} = \E \log_2(\D(j+r,1)) = \E \log_2{\Big(}\bigvee_{i=1}^{2^{j+r}} X(i) {\Big)}.
 $
 Therefore, 
 \begin{eqnarray}
  \E Y_{j+r} - (j+r)/\alpha - \E \log_2(\sigma_0 Z) &=& \E \log_2{\Big(}\frac{1}{2^{(j+r)/\alpha}}
 \bigvee_{i=1}^{2^{j+r}} X(i) {\Big)}
 - \E \log_2(\sigma_0 Z) \nonumber\\
 & = & \E \log_2(M_{n}) - \E \log_2(\sigma_0 Z),  \label{e:moments-1}
 \end{eqnarray}
 where $M_n:= n^{-1/\alpha}\bigvee_{i=1}^{n} X(i)$ and where $n:= 2^{(j+r)}$.  
 Corollary \ref{c:logs} implies that the right--hand side of \refeq{moments-1} is of order
 ${\cal O} (n^{-\beta/\alpha}) = {\cal O}(2^{-(j+r)\beta/\alpha}) = {\cal O}(2^{-r\beta/\alpha})$, 
 as $r\to\infty$, which in turn implies \refeq{moments}.

 {\it We now focus on proving \refeq{cov}.}  Let $i<j$ and recall that $N_{j+r} = [N/2^{j+r}],$ and 
 $N_{i+r} = [N/2^{i+r}].$ We also have that 
 \beq\label{e:cov-1} 
   \D(j+r,k) = \bigvee_{r = 1}^{2^{j-i}} \D(i+r,2^{j-i}(k-1)+i),\ \ \mbox{for all }k = 1,\ldots, N_{j+r}.
 \eeq
 Note that $2^{j-i}N_{j+r} \le N_{i+r}$ and therefore as in the proof of Proposition \ref{p:cov} above, we get
 \begin{eqnarray*}
  {\rm Cov}(Y_{j+r},Y_{i+r}) & = & \frac{1}{N_{j+r} N_{i+r}}\sum_{k_1=1}^{N_{j+r}}\sum_{k_2=1}^{N_{i+r}}
  {\rm Cov}(\log_2 \D(j+r,k_1),\log_2 \D(i+r,k_2)) \\
 &=& \frac{1}{N_{j+r} N_{i+r}}\sum_{k_1=1}^{N_{j+r}}\Sum_{\ell=1}^{2^{j-i}}
 {\rm Cov} {\Big(}\log_2 \D(j+r,k_1), \log_2 \D(i+r,2^{j-i}(k_1-1)+\ell){\Big)}.
 \end{eqnarray*}
 The second sum in the last expression involves only terms $\D(i+r,2^{j-i}(k_1-1)+\ell),$ for
 $\ell=1,\ldots,2^{j-i}$ since in view of \refeq{cov-1}, the independence of the $\D(i+r,k)$'s
 implies that ${\rm Cov}(\D(j+r,k_1),\D(i+r,k_2))=0$, for all $k_2$ outside the range 
 $2^{j-i}(k_1-1)+\ell,\ \ell=1,\ldots,2^{j-i}$.

 Now, by using the stationarity of the $\D(i+r,k)$'s and Relation \refeq{cov-1} again, we obtain from the 
 last relation that
 \begin{eqnarray}\label{e:cov-1.5}
  {\rm Cov}(Y_{j+r},Y_{i+r}) &=& \frac{2^{j-i}}{N_{i+r}} {\rm Cov}{\Big(} \log_2 {\Big(} \bigvee_{\ell=1}^{2^{j-i}}
  \D(i+r,\ell){\Big)}, \log_2 \D(i+r,1){\Big)} \nonumber\\
  &=& \frac{2^{j-i}}{N_{i+r}} {\rm Cov}{\Big(}\log_2{\Big(}M_n'\vee (2^{j-i}-1)^{1/\alpha} M_m''{\Big)},\log_2(M_n'){\Big)},
 \end{eqnarray}
 where $n:= 2^{i+r}$ and $m := (2^{j-i}-1)n$
 with  $M_n':= n^{-1/\alpha} \D(i+r,1) = n^{-1/\alpha}\bigvee_{\ell=1}^{n} X(\ell)$, and 
 $$
  M_m'' := m^{-1/\alpha} \bigvee_{\ell=2}^{2^{j-i}}\D(i+r,\ell) = m^{-1/\alpha} \bigvee_{\ell=n+1}^{n+m}X(\ell)
  \stackrel{d}{=} M_m'.
 $$
 Observe that the normalized maxima $M_n'$ and $M_m''$ are independent since they involve maxima of disjoint
 sets of $X(r)$'s.  Thus, by combining the results of Corollaries \ref{c:logs} and \ref{c:cov-logs}, we obtain
 that 
 \beq\label{e:cov-2}
  {\rm Cov}{\Big(}\log_2(M_n'\vee (2^{j-i}-1)^{1/\alpha} M_m''),\log_2(M_n'){\Big)} - \alpha^{-2}\psi(|i-j|)
  = {\cal O}{\Big(} 1/2^{r\beta/\alpha}{\Big)},\ \ r\to\infty,
 \eeq
 where $\psi$ is as in \refeq{psi}.
 Now, note that $N_{i+r} = 2^{j_2-i} N_{j_2+r} + q$, where $q <2^{j_2-i},\ q\in\bbN$.  This follows from the
 facts that $N_{i+r} = [N/2^{i+r}],$ $i=j_1,\ldots,j_2$ and $i\le j_2$.  Thus
 \beq\label{e:cov-3}
  \frac{N_{j_2+r}}{N_{i+r}} - 2^{i-j_2} = {\cal O}(1/N_{r}) = {\cal O}(2^r/N).
 \eeq
 Now, by applying Relations \refeq{cov-2} and \refeq{cov-3}, to \refeq{cov-1.5}, we obtain \refeq{cov}.  This
 completes the proof of the proposition.
 $\Box$
\end{proof}

\medskip
 The following theorem is the main result of the section.  It
 establishes the uniform convergence of the vector $Y_r$ to a
 normal vector and provides bounds on its rate of
 convergence.  The asymptotic normality of the estimators defined in
 \refeq{theta-GLS} is then an immediate consequence of this result (see Corollary \ref{c:an-alpha-sigma} below).
 
 \begin{theorem}\label{t:an} Suppose that the c.d.f.\  $F$ has the representation \refeq{F-sigma}
  and satisfies Conditions \ref{cond:C1} and \ref{cond:C2}, above.
  Let $\btheta = \{\theta_j\}_{j = j_1}^{j_2} \in \bbR^{m}\setminus\{{\mbox{{\bf 0}}}\},\
  m=j_2-j_1+1$ be an arbitrary fixed, non--zero vector and consider the linear
  combination $(\btheta,\bY_r) := \sum_{j=j_1}^{j_2} \theta_j Y_{j+r}$.

  Then, 
  \beq\label{e:an-Berry-Esseen}
   \sup_{x\in\bbR} {\Big|} \P\{ \sqrt{N_{j_2+r}}{\Big(}(\btheta,\bY_r) - (\btheta,\bmu_r){\Big)} \le x\} - \Phi(x/\sigma_\btheta) 
   {\Big|}
   \le  C_\btheta{\Big(}1/2^{r\beta/\alpha} + r 2^{r/2}/\sqrt{N}{\Big)},
  \eeq
  where $\Phi$ stands for the standard Normal c.d.f.\ and where $C_\btheta>0$ does not depend on $N$. Here
  $N_{j} = [N/2^{j}]$ denotes the number of coefficients $\D(j,k)$ available on scale $j$,
  $(\btheta,\bmu_r) :=  \sum_{j=j_1}^{j_2}\theta_j \mu_r(j)$ and 
  \beq\label{e:sigma-theta}
   \sigma_\btheta^2 = \alpha^{-2} (\btheta, \Sigma_1 \btheta) 
  := \alpha^{-2} \sum_{i,j=j_1}^{j_2} 
    \theta_i\Sigma_1(i,j)\theta_j > 0. 
  \eeq
\end{theorem} 
\begin{proof} Since $N_i = [N/2^i],\ i=1,\ldots,[\log_2 N]$, for all
 $j = j_1,\ldots,j_2$, and $r \in\bbN,\ r\le [\log_2 N] - j_2$, we have
 $N_{j+r} = 2^{j_2-j}N_{j_2+r} + q_j$, where $0\le q_j < 2^{j_2-j},\ q_j\in\bbN$.  
 Thus, for all $j=j_1,\ldots,j_2,$
\begin{eqnarray}
  \label{e:an-1}
  Y_{j+r} & = & \frac{1}{N_{j+r}} \sum_{k=1}^{N_{j_2+r}} \sum_{i=1}^{2^{j_2-j}} \log_2 \D(j+r, 2^{j_2-j}(k-1)+i)
  + \frac{1}{N_{j+r}} \sum_{i=1}^{q_j} \log_2 \D(j+r, 2^{j_2-j} N_{j_2+r} + i)\nonumber\\
  & =: & \frac{1}{N_{j_2+r}} \sum_{k=1}^{N_{j_2+r}} y_{j+r}(k) + R_j,
\end{eqnarray}
where $y_{j+r}(k) := N_{j_2+r} N_{j+r}^{-1} \sum_{i=1}^{2^{j_2-j}}
 \log_2 \D(j+r, 2^{j_2-j}(k-1)+i)$.  

\noi Therefore,
\beq\label{e:an-1.5}
  (\theta, Y_r) = \frac{1}{N_{j_2+r}} \sum_{k=1}^{N_{j_2+r}} \xi_r(k) + (\theta, R),
\eeq
where $\xi_r(k):= (\theta, y_r(k)),\ k=1,\ldots, N_{j_2+r},$ with 
$y_r(k) = \{y_{j+r}(k)\}_{j=j_1}^{j_2}$ and $R= \{R_j\}_{j=j_1}^{j_2}$.

Observe that the random vectors $y_r(k),\ k=1,\ldots, N_{j_2+r}$ are i.i.d.\   and
 independent from the remainder term $(\theta,R)$.  Indeed, this follows
from the fact that the $X(i)$'s are i.i.d.\  and because for any
$j=j_1,\ldots,j_2$, the random variable $y_{j+r}(k)$ depends only on the $X(i)$'s
with indices $2^{j_2+r}(k-1)+1\le i \le 2^{j_2+r}k$, $k=1,\ldots, N_{j_2+r}$, and
$R_j$ depends on the $X(i)$'s with indices $2^{j_2+r}N_{j_2+r}+1\le i \le N$. 

Thus, to prove \refeq{an-Berry-Esseen}, we proceed in two steps. First, we apply
the Central Limit Theorem to the first term on the right--hand side (r.h.s.) of \refeq{an-1.5}. 
Then, we will argue that the remainder term therein can be neglected.

\medskip
{\it Step 1.}  Note that the $\xi_r(k)$'s are i.i.d.\ but their distributions depend on $N$ and
 hence the ordinary C.L.T.\ does not apply.  The Berry--Esseen bound, however,
 (see e.g.\ Theorem V.2.4 in Petrov \cite{petrov:1995}) implies that
 \beq\label{e:an-2}
  \sup_{x\in\bbR} {\Big|} Q_{N,r}(x)  - \Phi(x) {\Big|} \le 
   A \frac{\E|\xi_r(1)-\E\xi_r(1)|^3}{\sigma_{\xi_r}^3} \frac{1}{\sqrt{N_{j_2+r}}},
 \eeq
 where
 $$
  Q_{N,r}(x) := \P {\Big\{}\frac{1}{\sigma_{\xi_r} \sqrt{N_{j_2+r}}} 
  \sum_{k=1}^{N_{j_2+r}}(\xi_r(k) - \E \xi_r(k)) \le x {\Big\}},
 $$
 $\Phi(x)$ denotes the standard Normal c.d.f., and where $A>0$ is an absolute constant. 
 This is so, provided that the variance
 $\sigma_{\xi_r}^2 := {\rm Var}(\xi_r(1))$ and the third moment $\E |\xi_r(1)|^3$
 of the $\xi_r(k)$'s are finite.
 
Observe first that, by \refeq{an-1.5} and by the independence of the $\xi_r(k)$'s from $R$,
\beq\label{e:an-sigma-xi}
  \sigma_{\xi_r}^2 = 
  N_{j_2+r}{\Big(} {\rm Var }( \theta,Y_r) - {\rm Var}(\theta,R) {\Big)} 
 = \sigma_\theta^2 + {\cal O}(1/2^{r\beta/\alpha}) + {\cal O}(2^r/N),   
\eeq
 where $\sigma_\theta$ is as in \refeq{sigma-theta}.
 Indeed, this follows from Proposition \ref{p:mom-cov} above, provided that
 ${\rm Var}(\theta,R)$ is negligible.  In view of \refeq{an-1}, however, since 
 $0\le q_j< 2^j\le 2^{j_2},\ j=j_1,\ldots,j_2$,
\begin{eqnarray}\label{e:an-3}
 {\rm Var}(\theta,R) &\le&  \frac{m^2 2^{j_2} }{N_{j_2+r}} \sum_{j=j_1}^{j_2}
 {\rm Var}(\log_2 \D(j+r,1))\nonumber\\
  &=& \frac{m^2 2^{j_2}}{N_{j_2+r}} \sum_{j=j_1}^{j_2}
 {\rm Var}(\log_2 (2^{-(j+r)/\alpha}\D(j+r,1)) ),
\end{eqnarray}
 where $m=j_2-j_1+1$. In the last relation, we used the inequality 
 ${\rm Var}(\eta_1+\cdots+\eta_m) \le m^2 ( {\rm Var}(\eta_1)+
 \cdots + {\rm Var}(\eta_m)),\ m\in\bbN$ and the fact that 
$$
 {\rm Var}( \log_2 \D(j+r,1) ) = {\rm Var}(\log_2 (2^{-(j+r)/\alpha}\D(j+r,1))).
$$
In view of \refeq{C-j}, however, by Corollary \ref{c:logs} below, 
the variances on the r.h.s.\ of \refeq{an-3} are bounded, as $r\to\infty$.
This implies that ${\rm Var}(\theta,R) = {\cal O}(2^r/N)$, which completes the proof
of \refeq{an-sigma-xi}.

\medskip
We now focus on bounding the term $\E|\xi_r(1)-\E\xi_r(1)|^3$ in \refeq{an-2}.
The inequality
 \beq\label{e:c-r}
  {\Big|}\sum_{i=1}^m x_i{\Big|}^p \le m^{0\vee(p-1)} \sum_{i=1}^m |x_i|^p,\ m\in\bbN,\ \ 
          \mbox{ valid for all } p,\ x_i\in\bbR,\ i=1,\ldots,m,
 \eeq
 implies
\begin{eqnarray}\label{e:an-3.5}
 \E |\xi_r(1) - \E\xi_r(1)|^3 &\le & 
 m^2 \sum_{j=j_1}^{j_2} |\theta_j|^3 \E |y_{j+r}(1)-\E y_{j+r}(1)|^3 \nonumber\\
 & \le & m^2 \sum_{j=j_1}^{j_2} |\theta_j|^3 \E {\Big|} \frac{1}{2^{j_2-j}} \sum_{i=1}^{2^{j_2-j}} 
    \log_2 \D(j+r,i) - \E \log_2 \D(j+r,1) {\Big|}^3\nonumber\\
 & \le & m^2 \sum_{j=j_1}^{j_2} \frac{|\theta_j|^3  }{2^{j_2-j}}
  \sum_{i=1}^{2^{j_2-j}}\E | \log_2 \D(j+r,i) - \E \log_2 \D(j+r,1)|^3\nonumber\\
 & = &  m^2 \sum_{j=j_1}^{j_2} |\theta_j|^3 \E | \log_2 \D(j+r,1) - \E \log_2 \D(j+r,1)|^3,
\end{eqnarray}
where $m=j_2-j_1+1$ and where the last bound follows from the Jensen's inequality. 
As in \refeq{an-3} above, we have that $\log_2 \D(j+r,1) - \E \log_2 \D(j+r,1)$ equals
$$ 
 \log_2(2^{-(j+r)/\alpha} \D(j+r,1)) - \E\log_2(2^{-(j+r)/\alpha} \D(j+r,1)),
$$
Therefore, by using inequality \refeq{c-r}, we get that the r.h.s.\ 
of \refeq{an-3} is bounded above by 
$$
 4m^2\sum_{j=j_1}^{j_2} |\theta_j|^3 {\Big(}\E |\log_2 (2^{-(j+r)/\alpha} \D(j+r,1))|^3
  + (\E |\log_2(2^{-(j+r)/\alpha} \D(j+r,1)|)^3{\Big)}.
$$
The last term is bounded, as $r\to\infty$, in view of \refeq{C-j} and Corollary
\ref{c:logs}.

\medskip
We have thus far shown that \refeq{an-2} holds with the r.h.s.\ being of
order ${\cal O}(1/\sqrt{N_{r}}),$ uniformly in $r$, that is,
\beq\label{e:an-4}
 \sup_{x\in\bbR} {\Big|} Q_{N,r}(x) - \Phi(x) {\Big|} \le C_\theta/\sqrt{N_r} = {\cal O}{\Big(}2^{r/2}/\sqrt{N} {\Big)}.
\eeq
We will now use this fact to prove \refeq{an-Berry-Esseen}.  

\medskip
{\it Step 2.} By \refeq{an-1.5}, the probability in \refeq{an-Berry-Esseen} equals
\beq\label{e:an-5}
 \E Q_{N,r} {\Big(}x/\sigma_{\xi_r} - \sqrt{N_{j_2+r}}((\theta,R) + \E \xi_r(1) 
 -(\theta,\mu_r))/\sigma_{\xi_r} {\Big)} =: \E Q_{N,r} {\Big(}x/\sigma_{\xi_r} - 
 \Delta_{N,r}{\Big)}.
\eeq
Indeed, this follows from the independence of the $\xi_r(k)$'s and the
remainder term $R$.

Now, by applying the triangle inequality, we obtain that the l.h.s.\ of
\refeq{an-Berry-Esseen} is bounded above by:
\begin{eqnarray}
 & & \sup_{x\in\bbR} \E {\Big|} Q_{N,r}(x/\sigma_{\xi_r} - \Delta_{N,r})
 - \Phi(x/\sigma_{\xi_r} - \Delta_{N,r}) {\Big|} + \sup_{x\in\bbR} \E {\Big|}\Phi(x/\sigma_{\xi_r} - \Delta_{N,r})
 - \Phi(x/\sigma_{\xi_r}) {\Big|}\nonumber\\
& & \ \ + \sup_{x\in\bbR} | \Phi(x/\sigma_{\xi_r}) - \Phi(x/\sigma_\theta)| =:
 A_1 + A_2 + A_3. \label{e:an-6}
\end{eqnarray}

In view of \refeq{an-4}, we have that
\beq\label{e:A1}
 A_1 \le \sup_{x\in\bbR} |Q_{N,r}(x) - \Phi(x)| =  {\cal O}{\Big(}2^{r/2}/\sqrt{N}{\Big)},
\eeq
as $N\to\infty$ and $N/2^r\to\infty$.

{\it Now, focus on the term $A_2$ in \refeq{an-6}.}  By using the mean value
theorem, for any $a<b,\ a,b\in\bbR$, we have that $|\Phi(a) - \Phi(b)| \le
|a-b|/\sqrt{2\pi}$.  Therefore (see \refeq{an-5}),
\beq\label{e:an-7}
A_2 \le \frac{1}{\sqrt{2\pi}} \E |\Delta_{N,r}| 
\le \frac{ \sqrt{N_{j_2+r}} }{\sqrt{2\pi}\sigma_{\xi_r}} {\Big(}
  \E |(\theta,R)| + \E|\xi_r(1) - (\theta,\mu_r)| {\Big)}. 
\eeq
As argued above, in view of \refeq{an-1}, we obtain by the triangle
inequality, that
\begin{eqnarray}
\E| (\theta,R)| &\le& \frac{\const}{N_{j_2+r}}
\Sum_{j=j_1}^{j_2} \E |\log_2 \D(j+r,1)| \nonumber\\
&\le& \frac{\const}{{N_{j_2+r}}}
\Sum_{j=j_1}^{j_2} \E |\log_2 (2^{-(j+r)/\alpha}\D(j+r,1))| + 
\const\frac{r}{{N_{j_2+r}}} = {\cal O}(r/{N_r}). \label{e:an-8}
\end{eqnarray}
The last relation follows by adding and subtracting the term $(j+r)/\alpha$,
and by applying Corollary \ref{c:logs} to the terms 
$\E |\log_2 (2^{-(j+r)/\alpha}\D(j+r,1))|$.

By \refeq{an-1}, $\E\xi_r(1) = \E (\theta, Y_r) - \E (\theta,R)$ and thus
by applying the triangle inequality, Proposition \ref{p:mom-cov} and Relation \refeq{an-8}, to
the second term in the r.h.s.\ of \refeq{an-7}, we obtain
\beq\label{e:A2}
 A_2 \le \const \sqrt{N_r}{\Big(} {r/N_r} + 1/2^{r\beta/\alpha} {\Big)} =  {\cal O}{\Big(} r2^{r/2}/N {\Big)} + 
 {\cal O}{\Big(} 1/2^{r\beta/\alpha}{\Big)}.
\eeq
Here, we also used the fact that $\sigma_{\xi_r} \to \sigma_\theta,\ \sigma_\theta>0,$
as $r\to\infty$ (see \refeq{an-sigma-xi} above).

{\it Consider now the term $A_3$ in \refeq{an-6}.}  As above, by using the mean value theorem, we obtain
\begin{eqnarray}\label{e:A3}
 A_3 &\le& \const |1/\sigma_{\theta} - 1/\sigma_{\xi_r}| = \const { |\sigma_\theta - \sigma_{\xi_r}| \over \sigma_{\theta}
 \sigma_{\xi_{r}}} \nonumber\\
  & =& {\cal O} ( 1/2^{r\beta/\alpha}) +  {\cal O} ( 2^{r}/N ),
\end{eqnarray}
as $r\to\infty$ and $N/2^r\to\infty$, where the last inequality follows form Relation \refeq{an-sigma-xi} above
and the fact that $\sigma_{\theta}^2 - \sigma_{\xi_r}^2 = (\sigma_{\theta} - \sigma_{\xi_r})
 (\sigma_{\theta} + \sigma_{\xi_r})$.

Now, by combining the bounds in Relations \refeq{an-6}, \refeq{A1}, \refeq{A2} and \refeq{A3}, we obtain
\refeq{an-Berry-Esseen}.  This completes the proof of the theorem.
$\Box$
\end{proof}

\medskip
Let now the scales $j_1\le j_2$ be fixed and let $r=r(N)\in\bbN,\ r+j_2\le [\log_2 N]$.   Theorem \ref{t:an}
shows that one can obtain consistent and asymptotically normal estimators of $H$ and $C = C(\sigma_0,\alpha)$,
as in the ideal Fr\'echet case \refeq{theta-GLS}.  Indeed, let $A = (a\ b)$ be as in \refeq{theta-GLS} and define $\what\theta_{\Sigma_1} =
(\what H_{\Sigma_1},\ \what C_{\Sigma_1})$
as in \refeq{theta-GLS} and $\alpha^{-2}\Sigma_1$ being the
asymptotic covariance matrix in Proposition \ref{p:mom-cov}.

By using \refeq{theta-GLS}, one can show that
\beq\label{e:H-C}
 \what H := \what H_{\Sigma_1} = \sum_{j=j_1}^{j_2} w_j Y_{j+r} \ \ \mbox{ and }\ \ 
 \what C := \what C_{\Sigma_1} = \sum_{j=j_1}^{j_2} v_j Y_{j+r} - r \what H_{\Sigma_1},
\eeq
where the $w_j$'s and the $v_j$'s are {\it fixed} weights such that
\beq\label{e:wj-vj}
 \Sum_{j=j_1}^{j_2} jw_j = \Sum_{j=j_1}^{j_2} v_j = 1\ \ \mbox{ and }\ \ 
 \Sum_{j=j_1}^{j_2} w_j = \Sum_{j=j_1}^{j_2} jv_j = 0.
\eeq
The following result establishes the asymptotic normality of these estimators.

\begin{proposition}\label{p:an-H-C}  Assume the conditions of Theorem \ref{t:an} hold.  
If $r=r(N)\in\bbN$ is
 such that $r2^{r}/N + 1/2^{r\beta/\alpha} \to 0,$ as $N\to\infty$, then for the estimators defined in \refeq{H-C},
 we have
 \beq\label{e:an-H-C}
   \sqrt{N_{j_2+r}} (\what H - H) \stackrel{d}{\longrightarrow} {\cal N}(0, H^2 c_w)
  \ \ \mbox{ and }\ \ 
   \sqrt{N_{j_2+r}/r} (\what C - C) \stackrel{d}{\longrightarrow} {\cal N}(0, H^2 c_w),
 \eeq
 as $N\to\infty$, where $c_w = \sum_{i,j=j_1}^{j_2} w_i w_j \Sigma_1(i,j)$ and where $C = C(\sigma_0,\alpha)$ is as in \refeq{the-C}.

\noi Moreover, 
$$
 \lim_{N\to\infty} N_{j_2+r}{\rm Var}(\what H) = \lim_{N\to\infty} r^{-1}N_{j_2+r} {\rm Var}(\what C) =  H^2 c_w.
$$
\end{proposition}
\begin{proof}  The first convergence in \refeq{an-H-C} follows directly from Theorem \ref{t:an} by setting
$\theta_j := w_j,\ j=j_1,\ldots,j_2$.  Indeed, since $\mu_r(j)= (j+r)/\alpha + C$, Relation \refeq{wj-vj} implies
that 
$$
(\theta,\mu_r) = \Sum_{j=j_1}^{j_2} w_j ((j+r)/\alpha + C) = 1/\alpha \equiv H.
$$
Thus, for $\what H = (\theta,Y_r) = \sum_{j=j_1}^{j_2} w_j Y_{j+r}$, by Relation \refeq{an-Berry-Esseen}, we obtain that
$$
 \sup_{x\in\bbR} | \P\{ \sqrt{N_{j_2+r}}(\what H - H) \le x\} - \Phi(x/\sigma_w)| \longrightarrow 0,
$$
as $N\to\infty$.  This implies the asymptotic normality of $\what H$ in \refeq{an-H-C},
where in view of \refeq{sigma-theta} $\sigma_w^2 = H^2(w,\Sigma_1 w) = H^2 \sum_{i,j=j_1}^{j_2} w_iw_j \Sigma_1(i,j)$.

{\it We now focus on the estimator $\what C$.}  By setting $\theta_j:= v_j,\ j=j_1,\ldots,j_2$, we get by using \refeq{wj-vj}
that
$$
(\theta,\mu_r) = \sum_{j=j_1}^{j_2}((j+r)/\alpha +C) v_j = r/\alpha + C.
$$ 
On the other hand, in view of \refeq{H-C},
$$
(\theta,Y_r) = \Sum_{j=j_1}^{j_2} v_j Y_{j+r} = \what C + r \what H
$$
and thus
\beq\label{e:p:an-H-C-1}
 \what C - C = (\theta,Y_r) - (\theta,\mu_r) - r(\what H - H). 
\eeq
We have already shown that the term $(\what H -H) $ above is asymptotically normal and by Theorem \ref{t:an}
the term $(\theta,Y_r) - (\theta,\mu_r)$ in \refeq{p:an-H-C-1} is also asymptotically normal.  Since $r=r(N)\to \infty$,
the second term in the r.h.s.\ of \refeq{p:an-H-C-1} dominates in the limit.  
This implies that second convergence in \refeq{an-H-C}.

To complete the proof, observe that by Proposition \ref{p:mom-cov}, $N_{j_2+r} {\rm Var}(\what H) \to \sigma_w^2 = H^2 c_w,$
as $N\to\infty$.  We now consider the variance of $\what C -C$ in \refeq{p:an-H-C-1}, and apply the inequality
$$
 {\rm Var}(\xi)- 2({\rm Var}(\xi){\rm Var}(\eta))^{1/2} +  {\rm Var}(\eta) \le {\rm Var}(\xi - \eta)
\le {\rm Var}(\xi) + 2({\rm Var}(\xi){\rm Var}(\eta))^{1/2} +  {\rm Var}(\eta)
$$
with $\xi:= (\theta,Y_r) - (\theta,\mu_r)$ and $\eta:= r(\what H - H)$.  
Since ${\rm Var}(\eta)$ dominates ${\rm Var}(\xi)$, in the limit, we obtain that
$r^{-1}N_{j_2+r} {\rm Var}(\what C) \to \sigma_w^2 = H^2 c_w$, as $N\to\infty$.
$\Box$
\end{proof}

\begin{corollary}\label{c:an-alpha-sigma} Assume the conditions of Theorem \ref{t:an} hold. 
Define the estimators 
$$
  \what \alpha := 1/\what H\ \ \mbox{ and }\ \ \what \sigma_0 := 2^{\what C - (\E \log_2 Z)/\what \alpha},
$$
where $Z$ is a $1-$Fr\'echet random variable with unit scale coefficient.  Then with $r=r(N)$ as in Proposition
\ref{p:an-H-C}, we have
 \beq\label{e:an-alpha-sigma}
   \sqrt{N_{j_2+r}} (\what \alpha - \alpha) \stackrel{d}{\longrightarrow} {\cal N}(0, \alpha^2 c_w)
  \ \ \mbox{ and }\ \ 
   \sqrt{N_{j_2+r}/r} (\what \sigma_0 - \sigma_0) \stackrel{d}{\longrightarrow} {\cal N}(0, (\ln 2)^2 \sigma_0^2 \alpha^{-2} c_w).
 \eeq
\end{corollary}

\medskip
\noi This result follows from Proposition \ref{p:an-H-C} by an application of the Delta-method.

\medskip
Most heavy--tailed  distributions used in applications
satisfy Condition \ref{cond:C1}, but some do not satisfy Condition \ref{cond:C2}.
Indeed, \refeq{C2} implies that $\E |X|^p 1_{\{X\le 1\}} <\infty$, for all $p\in\bbR$, which is rather stringent.
Nevertheless, the results of Proposition \ref{p:an-H-C} and Corollary \ref{c:an-alpha-sigma} continue to hold even if
Condition \ref{cond:C2} is not satisfied and even if the $X(i)$'s can take negative values.  This is so, because block--maxima
become strictly positive as the block--size grows.  We make this more precise in Proposition \ref{p:coupling} below.

 Now, for convenience, introduce a special value $*$ and suppose that our statistics take values in the extended
real line $\bbR^*:= \bbR \cup \{*\}$.  If a statistic is not well--defined (because it involves $\log_2 x$ for $x\le 0$, for example),
we assign to it the special value $*$.   The set $\{*\} \subset \bbR^*$ is considered as both closed and open in the topology
of $\bbR^*$ and the topology of $\bbR\subset \bbR^*$ is the same as that of the real line.   Therefore, the statistics $Y_j$ in
\refeq{Yj} and the estimators $\what H$ and $\what C$ in \refeq{H-C}, become proper random variables 
which can sometimes take the value $*$ if some of the $X(i)$'s are negative.  

The following result shows that, asymptotically, the estimators $\what H$ and $\what C$ become real--valued with
probability one, provided that $\ln(N)/2^{r(N)}\to 0$, as $N\to\infty$.  

\begin{proposition}\label{p:coupling}  Suppose that the c.d.f.\ $F$ has the representation \refeq{F-sigma} and satisfies Condition \ref{cond:C1},
where $F(0)$ is not necessarily zero.  Let also $r=r(N)\in\bbN$, $\what H$ and $\what C$ be as in \refeq{H-C}.  If
$ \ln(N)/2^{r(N)} \longrightarrow 0,\ N\to\infty$, then
\beq\label{e:p:coupling}
 \P (\{\what H = *\}) + \P(\{\what C = *\}) \longrightarrow 0,\ \ \mbox{ as } N \to\infty.
\eeq
If in addition $r2^r/N + 1/2^{r\beta/\alpha} \to 0$, as $N\to\infty$, then the convergences \refeq{an-H-C} 
and \refeq{an-alpha-sigma} continue to hold.
\end{proposition} 
\begin{proof} Let $X(i),\ i\in\bbN$ be i.i.d.\ with c.d.f.\ $F$ and let $x_0>0$ be arbitrary.
Define the truncated variables $\wtilde X(i) := X(i) 1_{\{X(i)> x_0\}} + x_0 1_{\{X(i)\le x_0\}},\ i\in\bbN$ and observe that they
are i.i.d.\ with c.d.f.\ $\wtilde F(x) :=  F(x),\ x\ge x_0$ and $\wtilde F(x) = 0,\ x < x_0$.
Thus, $\wtilde F(x)$ has a representation as in \refeq{F-sigma} with the 
function $\sigma^\alpha(x)$ replaced by
$$ \wtilde \sigma^\alpha(x) = \infty 1_{(-\infty, x_0)}(x) + \sigma^\alpha(x) 1_{[x_0,\infty)}(x),$$
where $\sigma^\alpha(x)$ is the function involved in the corresponding representation of $F(x)$.

Consider the statistics $\wtilde \D(j,k)$ and $\wtilde Y_j$ defined as in \refeq{C-j} and \refeq{Yj} with
$X(i)$'s replaced by $\wtilde X(i)$'s.  Let also $\wtilde H$ and $\wtilde C$ be the corresponding statistics
defined as in \refeq{H-C} with $Y_j$'s replaced by $\wtilde Y_j$'s.
Observe that $\wtilde F$ satisfies Condition \ref{cond:C1} and also trivially Condition \ref{cond:C2} since $x_0>0$ and
$\wtilde \sigma^\alpha(x) = \infty$ for all $x\in (0,x_0)$.  Therefore, the results of Proposition \ref{p:an-H-C} apply to the
statistics $\wtilde H$ and $\wtilde C$.  We will now show that the statistics $\what H $ and $\what C$, which may not
be always real--valued random variables (i.e.\ can take the special value $*$) coincide with the statistics
$\wtilde H$ and $\wtilde C$, eventually.

Let $1\le j_0\le \log_2 N,\ j\in \bbN$.  Observe that  the event 
$${\cal C}_{j_0}:= \{ \wtilde \D(j_0,k) = \D(j_0,k),\ k=1,\ldots,N_{j_0}\}$$
implies the events ${\cal C}_{j} = \{ \wtilde \D(j,k) = \D(j,k),\ k=1,\ldots,N_{j}\}$, for all $j_0 \le j \le \log_2 N$ and 
in particular the events $\{\wtilde Y_j = Y_j\},\ j\ge j_0$.  Thus, the statistics $\wtilde H$ and $\what H$
(and $\wtilde C$ and $\what C$, respectively) coincide on the event ${\cal C}_{j_1+r}$.  Thus, to complete the 
proof of the proposition, it is sufficient to show that $\P( {\cal C}_{j_1+r})\to 1$, as $N\to\infty$.

Let $j_0:=j_1+r$ and observe that by independence, 
$$
\P({\cal C}_{j_0}) = \P\{ \wtilde \D(j_0,1) = \D(j_0,1)\}^{N_{j_0}} = {\Big(}1 - F(x_0)^{2^{j_0}}  {\Big)}^{N_{j_0}}.
$$
In view of Condition \ref{cond:C1}, $p_0:= F(x_0)<1$ and hence
$$
\ln \P({\cal C}_{j_0}) = N_{j_0} \ln(1- p_0^{2^j}) =  - \frac{N}{2^{j_0}} p_0^{2^{j_0}}(1+o(1)),\ \mbox{ as }j_0 \to\infty.
$$
Since $p_0<1$, the first convergence in \refeq{p:coupling} implies that $N p_0^{2^{j_1+r(N)}} \to 0,$ as $N\to\infty$,
and hence $\P({\cal C}_{j_1+r(N)})\to 1,$ as $N\to\infty$.  We have thus shown that \refeq{an-H-C} holds.  
Relation \refeq{an-alpha-sigma} follows from \refeq{an-H-C} by using the Delta--method.
$\Box$
\end{proof}

\medskip
\noi{\bf Remarks:}
\begin{enumerate}
 \item Observe that in view of \refeq{theta-GLS}, $\what H_{\Sigma_1} = \what H_{\phi \Sigma_1}$ and
       $\what C_{\Sigma_1} = \what C_{\phi \Sigma_1}$, for any $\phi>0$.  That is, one can compute, in practice,
       the generalized least squares estimators $\what H$ and $\what C$ without having to use a plug--in estimator
       for $\alpha$ in \refeq{cov} (see also the Remarks in Section \ref{s:GLS}).
 
 \item The constants $c_w$ appearing in Proposition \ref{p:an-H-C} and Corollary \ref{c:an-alpha-sigma}
       are given in Table \ref{tab:cw} below.  We now comment on the optimal rate in these asymptotic results.
      
       Proposition \ref{p:mom-cov} indicates that the bias of the estimator $\what H$ in \refeq{an-H-C} is of order
       ${\cal O}(1/2^{r\beta/\alpha})$. On the other hand, the standard error of $\what H$ is of order ${\cal O}(2^r/N)$.
       By balancing these orders, we obtain that 
       $$2^r = 2^{r(N)} \propto N^{\alpha/(2\beta+\alpha)} $$
       yields the {\it optimal order} of the mean squared error (m.s.e.)
       $\E (\what H - H)^2$, and a corresponding rate of convergence
       $$
        2^{r/2}/\sqrt{N} = {\cal O}(1/N^{\beta/(2\alpha+\beta)})
       $$
       to the limit distribution of $\what H$ in \refeq{an-H-C}.  

       Hall \cite{hall:1982} (see Theorem 2 therein) obtained the same optimal
       order of convergence for the Hill--type estimators under the
       following semi--parametric assumptions on the tail of $F$:
       \beq\label{e:hall-F}
           1- F(x) = c_1 x^{-\alpha}(1+ c_2 x^{-\beta} + o(x^{-\beta})),\ \ \mbox{ as } x\to\infty,\ \alpha, \beta>0.
       \eeq
       A Taylor expansion shows that this tail behavior corresponds to Condition \ref{cond:C1} above in
       the case when $0<\beta\le \alpha$.
       {\it Note that in Hall \cite{hall:1982} the parameter $r$ corresponds to $N/2^r$ in our case.}

       Observe that Theorems 1 and 2 in Hall \cite{hall:1982} involve also asymptotic normality results for
       the scale parameter $c_1$ in \refeq{hall-F}.  These results are similar to those about $\what C$ in 
       Proposition \ref{p:an-H-C}.  Note in particular the presence of the logarithmic in $N$ factor $r=r(N).$

\item  The optimal rate in the previous remark may not be improved, in general.  Indeed, by 
       Proposition \ref{p:rate-exact} the rate of the bias is exact if $\sigma^\alpha(x) - \sigma_0^\alpha \sim c_1 x^{-\beta},\ 
       x\to\infty,\ c_1\not=0$.   This is typically the case in practice (see the Examples above).
       Relation \refeq{cov} also implies that the order of the variance of $\what H$ is precisely
       ${\cal O}(1/\sqrt{N_r})$, and cannot be improved.  

       Furthermore, the rate in the Berry--Esseen bound may not be improved, in general (see e.g.\ Ch.\ V.2 in Petrov
       \cite{petrov:1995}).  Thus,  the result of Theorem \ref{t:an} is optimal in our setting.
       
\item  Consider the case of optimal m.s.e.\  of $\what H$, that is, $2^{r} \propto N^{\alpha/(2\beta+\alpha)}$.
       Observe that the r.h.s.\ in \refeq{an-Berry-Esseen} is up to the logarithmic in $N$ factor of $r(N)$ of the same
       order as the root--m.s.e.\ $(\E (\what H - H)^2)^{1/2}$.  This indicates that the precision (in terms of coverage probability)
       of the confidence intervals for $H$ based on the asymptotic distribution for $\what H$ will be of order at least
       ${\cal O}(1/N^{\beta'/(2\alpha+\beta)})$ for any $\beta'\in(0,\beta)$.

 \item Even though the estimators $\what \alpha$ and $\what \sigma_0$ in Corollary \ref{c:an-alpha-sigma} are asymptotically
 normal, it is not a good idea to use their asymptotic distributions to construct confidence intervals for $\alpha$ and $\sigma_0$.
 Indeed, for simplicity consider the ideal Fr\'echet case.  In this case, the estimator $\what H$ is unbiased and hence the
 estimator $\what \alpha = 1/\what H$ is {\it biased}.  Moreover, since the variance of the random variable $1/X$, where
 $X$ has Normal distribution is infinite, we expect that ${\rm Var}(\what \alpha)$ does not converge to the asymptotic
 variance of $\what \alpha$ in \refeq{an-H-C}.  In our experience, the distribution of $\what \alpha$ tends to be skewed in practice.
 Therefore, one can get better confidence interval estimates for $\alpha$ by using {\it inversion} from the corresponding
 confidence intervals for $H$.  For example, $((\what H + z_{p} \what H \sqrt{c_w}/\sqrt{N_{j_2+r}})^{-1}, 
 (\what H - z_p \what H \sqrt{c_w}/\sqrt{N_{j_2+r}})^{-1})$ is an asymptotically correct $100(1-p)\%$ confidence
 interval for $\alpha$, where $z_p:= \Phi^{-1}(1-p/2),\ p\in(0,1)$.  As indicated in the previous remark the error in the coverage
 probability of this interval is of order ${\cal O}(1/N^{\beta'/(2\alpha+\beta)})$ for any $\beta'\in(0,\beta)$, if m.s.e.--optimal
 $r$'s are chosen.

\end{enumerate}

\section{Performance evaluation and data analysis}
 \label{s:performance}

 \subsection{Typical models: small and large sample properties}
 \label{s:performance:typical}

We study the performance of the max self--similarity estimators when the data
are heavy--tailed but deviate from the ideal  Fr\'echet case. Specifically, 
given a sample of size $N=2^{n},\ n\in\bbN$, 
the GLS estimators $\what H = \what H(j_1,j_2)$ and $\what \alpha = \what \alpha (j_1,j_2) = 1/\what H$
are computed for a range of scales $j_1\le j \le j_2.$  We choose here $j_2=n$ as the maximal 
 available scale and focus on optimal $j_1$'s in the sense of mean squared error.  Namely, 
 we let  \beq\label{e:j1-opt}
  j_1^{opt} := \mathop{{\rm Argmin}}_{j_1,\, 1\le j_1\le j_2} \E(\what H(j_1,j_2) - H)^2,
 \eeq
 where the last expectation is computed from samples of independent
 realizations of the estimators $\what H$.

\begin{figure}[th!]
\begin{center}
\includegraphics[width=4.5in]{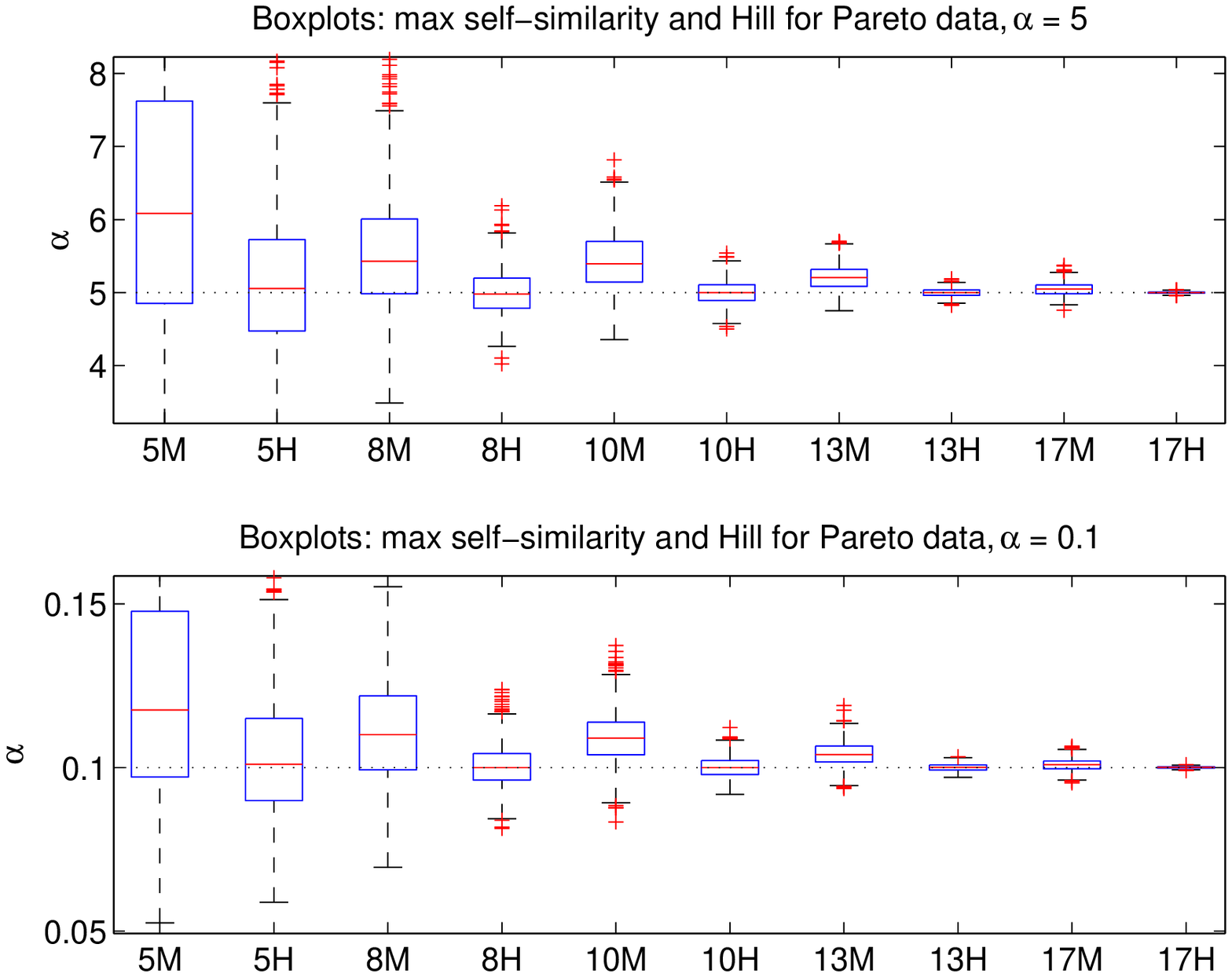}
{\caption{\label{fig:pareto_hill_vs_max_ss} \small 
Boxplots of $1,000$ independent realizations of max self--similarity and Hill estimators
for different sample sizes from Pareto distributions with $\alpha=5$ ({\it top panel}) and $\alpha=0.1$
 ({\it bottom panel}) are shown. The labels $nM$ and $nH$ correspond to sample size $2^n$ of max 
self--similarity and Hill estimators, respectively.
The Hill estimators were computed by using \refeq{alpha-k-Y}
with $k=2^n-1$, and the max self--similarity estimators are based on a range of scales 
$j_1\le j\le j_2=n$, where $j_1$ was chosen to minimize the mean squared error. }}
\end{center}
\end{figure}

 We first compare the max self--similarity estimators to the classical Hill estimator over
 Pareto data with unit scale, i.e.\ with c.d.f.\ $F(x) = 1-x^{-\alpha},\ x\ge 1$.
 In this case, the Hill estimator corresponds to the maximum likelihood estimator.
 Figure \ref{fig:pareto_hill_vs_max_ss} indicates that, as expected, the Hill estimator outperforms
 the max self--similarity estimator.  However, as seen from the box--plots, the max self--similarity
 estimator works relatively well for small, moderate and large samples and essentially keeps up with 
 the Hill estimators.  In fact, as the sample size grows the max self--similarity estimator improves 
 almost at the same rate as the Hill estimator.  Here the max self--similarity estimator was computed
 by using the range of scales $j_1^{opt}\le j \le j_2,$ where $j_2 = \log_2 N$ and $j_1^{opt}$ is
 as in \refeq{j1-opt}.

\begin{figure}[th!]
 \begin{center}
 \includegraphics[width=4.5in]{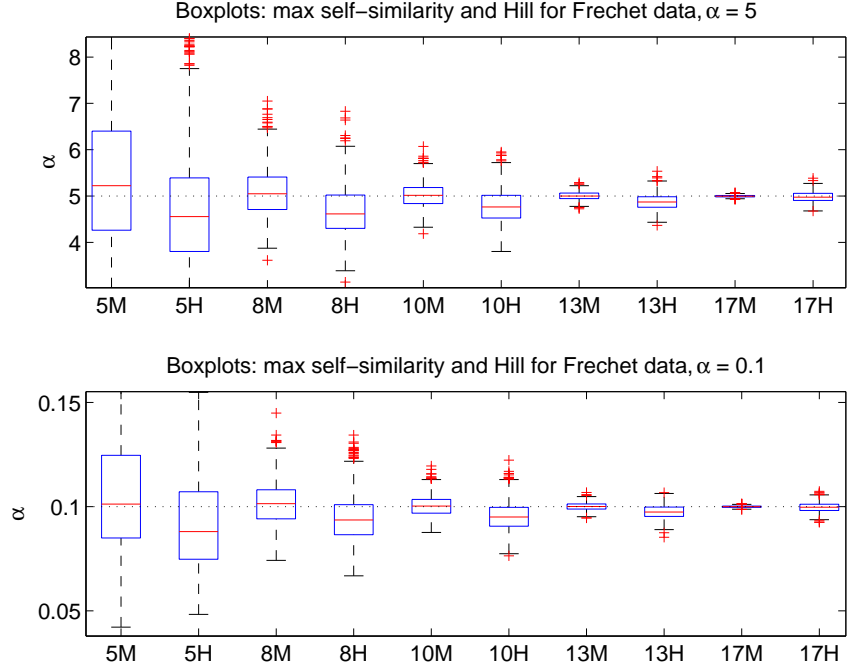}
  {\caption{\label{fig:frechet_hill_vs_max_ss} \small 
  Boxplots of $1,000$ independent realizations of max self--similarity and Hill estimators
  for different sample sizes from Fr\'echet distributions with $\alpha=5$ ({\it top panel}) and $\alpha=0.1$
  ({\it bottom panel}) are shown. The labels $nM$ and $nH$ correspond to sample size $2^n$ of max 
  self--similarity and Hill estimators, respectively.  
  The Hill estimator were computed by using an optimal value for $k$ in
  \refeq{alpha-k-Y}, which yields the smallest mean squared error.  The max self--similarity estimators
  were computed from the entire range of scales $j$.}}
 \end{center}
\end{figure}  

 In Figure \ref{fig:frechet_hill_vs_max_ss}, we compare the performance of the max self--similarity and
 the Hill estimators for Fr\'echet data.  The parameter $k$ in \refeq{alpha-k-Y} of the Hill estimator
 was chosen to minimize the mean squared error of the statistics $1/\what \alpha_H(k)$, by analogy 
 with \refeq{j1-opt}.  Now, the entire range of scales $j_1=1\le j_2=\log_2 N$ was used to compute
 the max self--similarity estimators.  Observe that as compared to the case of Pareto data (see Figure
 \ref{fig:pareto_hill_vs_max_ss}), now the roles of the two estimators are reversed.  As expected,
 the max self--similarity estimator works best in the Fr\'echet setting and dominates the Hill estimator.
 In fact, the method of choosing the parameter $k$ here is unusually favorable to the Hill estimator since
 it is not based on examining and determining a range where the Hill plot is constant.  It is well known that 
 in practice, the Hill plot is quite volatile and the resulting choice of $k$ based on this
 plot would yield far more biased estimators than the ones shown in Figure \ref{fig:frechet_hill_vs_max_ss}.

 \medskip
 We now examine the max self--similarity estimators in more detail 
 when the data are drawn from a stable and a $t-$distribution.
 Tables \ref{tab:stable} and \ref{tab:t_dist} below, indicate that the
 estimators $\what H_{opt}:=\what H(j_1^{opt},j_2)$ work well in
 practice for a variety of sample sizes and parameter values.
 Their performance is particularly good in the stable context.
 The performance in the case of $t-$distributions is comparable with the
 stable cases when the heavy--tail exponent $\alpha$ is not large. Notice that
 $\alpha$ corresponds to the degrees of freedom of the $t-$distribution and therefore
 as $\alpha$ grows, the $t-$distribution gets
 closer to the Normal distribution. Although it it still heavy tailed, most of the body
 of the distribution is not and therefore the quality of the tail estimators deteriorates.

 Table \ref{tab:stable} indicates that the max self-similarity estimator outperforms
 the Hill estimator for stable distributions with $\alpha \le 1$ and that the two estimators
 are comparable for $1<\alpha<2$.  The Hill estimator
 is slightly better than or comparable to the max self-similarity one
 for the t-distributions with low $\alpha$'s and slightly worse or comparable
 for moderate and large $\alpha$'s (Table \ref{tab:t_dist}).  

 The MSE--optimal choice of the parameter $k$ is unrealistically favorable
 to the Hill estimator.  In practice, these choices of $k$ typically do not 
 correspond to constant regions in the Hill plot.  On the other hand
 the MSE-optimal values of $j_1$ usually correspond to the knee in the 
 max--spectrum plot, which can be identified in practice (either visually
 or automatically).  These observations suggest that in reality the max 
 self--similarity estimators are more reliable and accurate than estimators
 based on the Hill plot.
 
\subsection{On the selection of the scales $j_1$ and $j_2$}
 \label{s:cut-off}

 In the ideal case of $\alpha$-Frechet data, the max--spectrum plot of $Y_j$ is almost
perfectly linear in $j$ (see Figure \ref{fig:max_spectrum}).  However, most real
data sets deviate from the ideal case and thus the max--spectrum 
becomes linear only over a range of relatively large scales $j$.  The selection of
an appropriate range of scales $j_1\le j\le j_2$, where the max self--similarity estimators are
computed, becomes an important practical problem. Because of \refeq{Yj-sim}, one
can always choose $j_2 = [\log_2 N]$ to be the largest available scale and the scale
$j_1$ can be chosen by visual inspection, a strategy that work fairly well in practice.
Nevertheless, we also propose an {\it automatic procedure} for choosing
the scale $j_1$, which turns out to also work well in practice.  It relies on the following 
simplifying assumptions:

\smallskip
\noi{\bf Assumption 1.}  The vector $Y_j,\ j=1,\ldots,j_2$ follows a {\it multivariate 
Normal distribution}.  

\smallskip
\noi{\bf Assumption 2.}  The covariance matrix $\Sigma_\alpha(1,j_2;N)= \alpha^{-1}\Sigma_1(1,j_2;N)$
of the vector $Y = \{Y_j\}_{j=1}^{j_2}$ is given by \refeq{p:cov}.

\medskip
These assumptions are valid asymptotically, provided that $N_{j_2}\to\infty$ (Theorem \ref{t:an} and Proposition \ref{p:mom-cov}).
Since the $N_{j}$'s grow exponentially fast as $j$ decreases, choosing $j_2$ 
as the largest available scale $[\log_2 N]$ is not critical in practice.
Let now $\what H(j_1,j_2)$ denote the GLS estimate of $H=1/\alpha$, computed over the range of scales $j_1\le j \le j_2$ as in \refeq{theta-GLS} (see also \refeq{H-C}).

\begin{figure}[t!]
\begin{center}
\includegraphics[height=2in,width=2.5in]{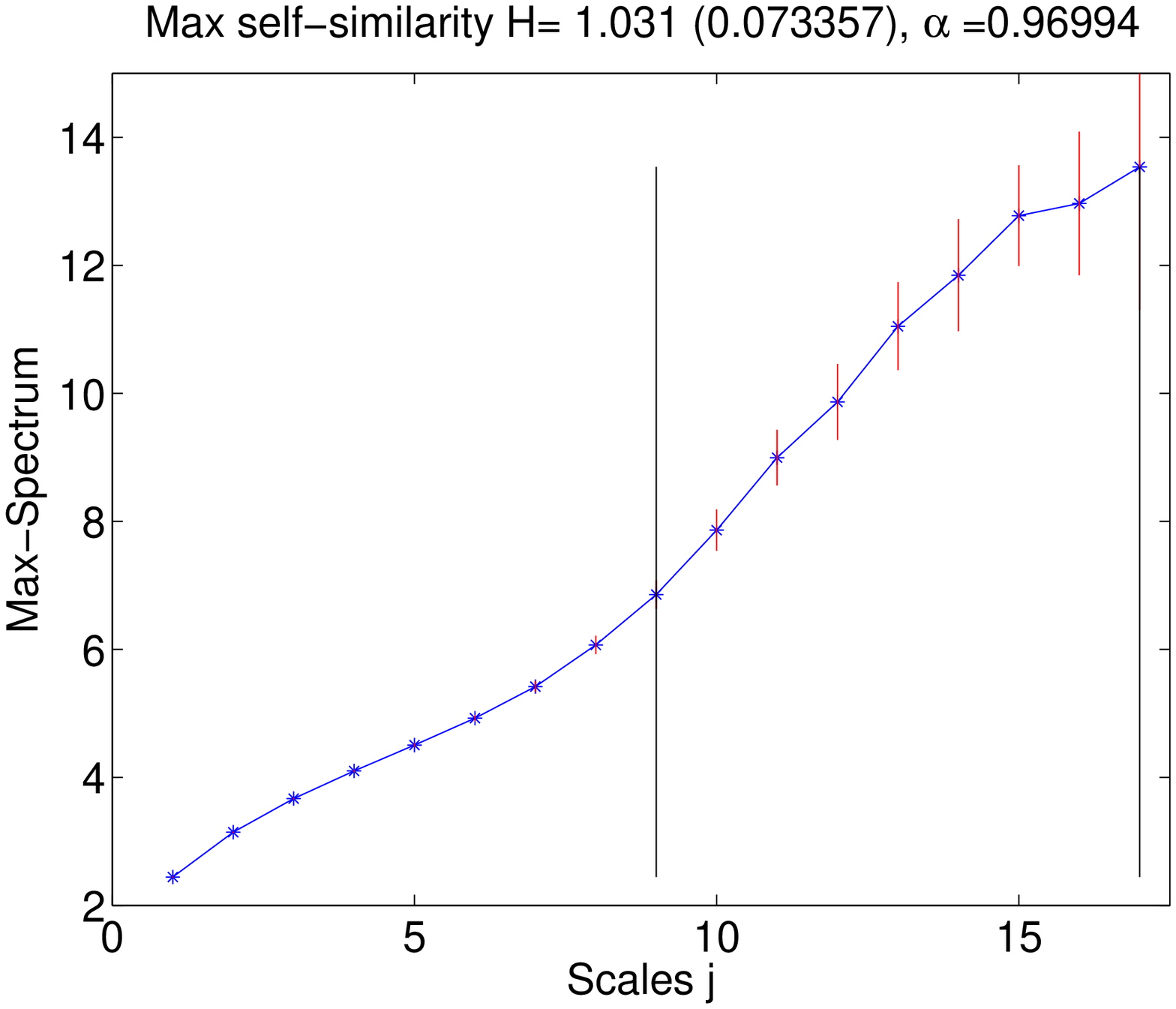} \hskip .2 in
\includegraphics[height=2in,width=2.5in]{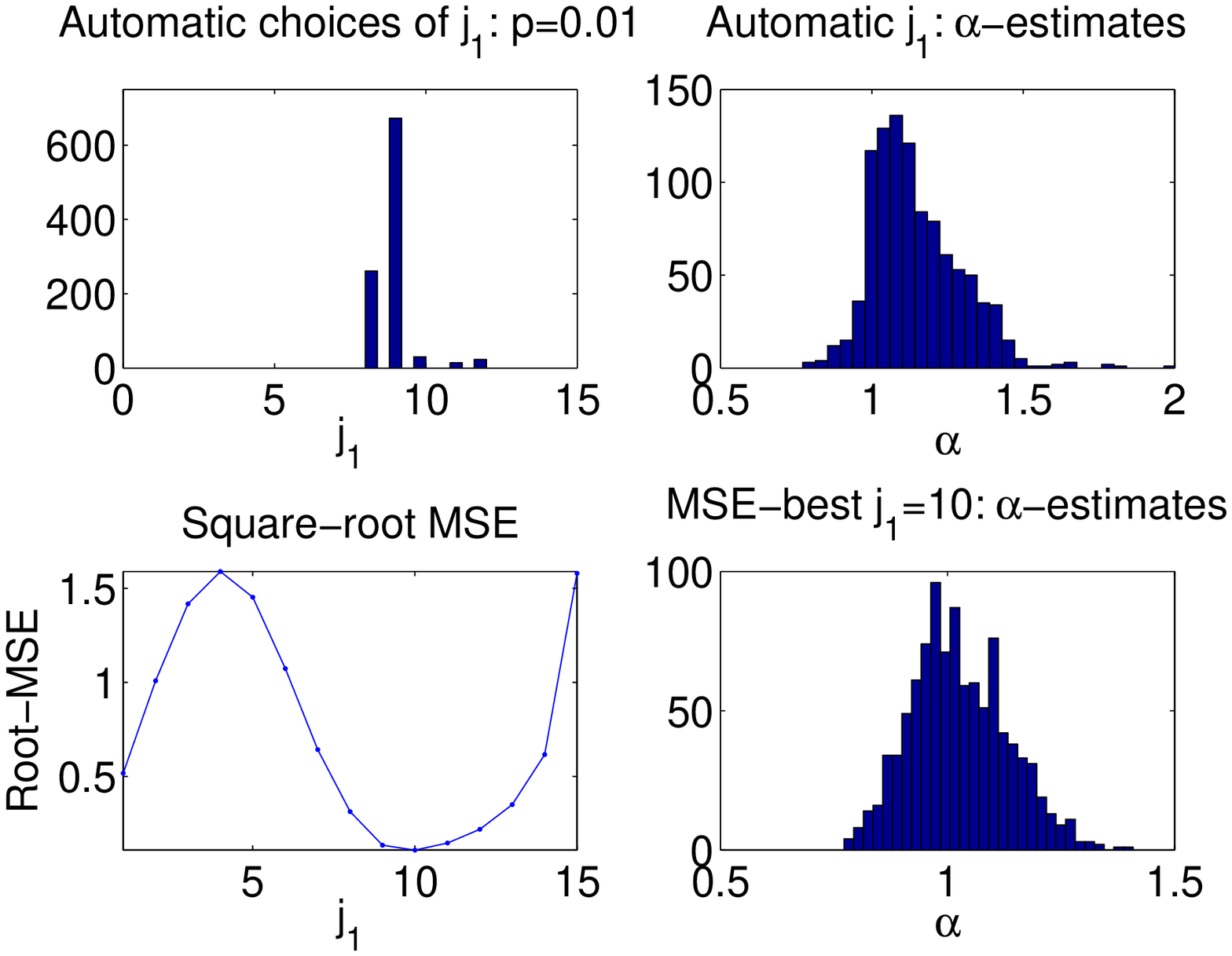}
{\caption{\label{fig:cut-off} \small  Mixtures of $\alpha-$Fr\'echet ($10\%$) and Exponential of mean $5$ ($90\%$)
were simulated. The heavy--tail exponent is $\alpha=1$ and the sample sizes are 
$N=2^{17} = 131,072$.  
{\it Left panel:} max--spectrum of a typical sample. {\it Right panel:}  $1,000$ independent replications of 
the GLS max self--similarity estimators were obtained, where automatic selection for the parameter $j_1$ was used
with $p=0.01$ and $b=4$.  The top--right graph shows a histogram of the resulting selections of $j_1$.  The bottom--right
graph shows the root--mean squared error of the estimators $\what H = 1/\what\alpha$.  The top--left and top--right plots 
shows histograms of the $\what \alpha$ estimates obtained by using automatically selected $j_1$'s and with
 $j_1=j_1^{opt}=10$, respectively.
 }}
\end{center}
\end{figure}

\medskip
\noi\emph{\bf Algorithm}

\medskip
\noi
{\it Tunning parameters:}

\noi{\hskip 0.5in}\parbox[t]{5.7in}{Pick a relatively small {\it significance threshold} $p \in (0,1)$ 
(e.g.\ $p=0.1$ or $0.01$) and an integer $b$ called {\it back--start parameter} (e.g.\ $b=3$ or $4$ for moderate sample sizes).  Set $j_2 := [\log_2 N]$ and $j_1:= \max\{1, j_2-b\}$.}

\smallskip
\noi
{\it Step 1. }\parbox[t]{5.7in}{If $j_1=1$ then stop, 
                else calculate $\what H_{\rm new} = \what H(j_1-1,j_2)$ and $\what H_{\rm old} = \what H(j_1,j_2)$.}

\smallskip
\noi
{\it Step 2. }\parbox[t]{5.7in}{Let $w_{\rm new}$ and $w_{\rm old}$ be vectors of weights as in \refeq{H-C}, such that
$\what H_{\rm new} = (w_{\rm new},Y)$ and $\what H_{\rm old} = (w_{\rm old},Y)$, where $Y = \{Y_{j}\}_{j=1}^{j_2} 
\in \bbR^{j_2}$ and where the vectors $w_{\rm new},\ w_{\rm old} \in\bbR^{j_2}$ are appropriately padded with zeros.   
Consider the quantity:
$$
S_1:= {\Big(}(w_{\rm new}-w_{\rm old}), \Sigma_1(1,j_2;N) (w_{\rm new}-w_{\rm old}) {\Big)}^{1/2}
$$ 
Now, consider the approximate $(1-p)-$level confidence interval for $\E(\what H_{new} - \what H_{old})$:
$$
 {\Big(} \what H_{new} - \what H_{old} - z_{p/2}  \what H_{\rm old} S_1, 
         \what H_{new} - \what H_{old} + z_{p/2} \what H_{\rm old} S_1 {\Big )},
$$
where $z_{p/2} = \Phi^{-1}(1-p/2)$ is a $(1-p/2)-$th quantile of the standard Normal distribution.}

\medskip
\noi
{\it Step 3. }\parbox[t]{5.7in}{ If zero is contained in the confidence interval computed in {\it Step 2}, then set 
$j_1 := j_1 -1$ and go to {\it Step 1} otherwise stop and report the selected $j_1$ and $\hat{\alpha}:= 1/\what H_{\rm old}$.}\\

\medskip
The choice of tunning parameters $p$ and $b$ and the validity of the above simplifying assumptions is addressed 
in Stoev et al. \cite{stoev:michailidis:hamidieh:taqqu:2006P}.  In Figure
\ref{fig:cut-off},  we briefly demonstrate the performance of the above automatic selection procedure for
a mixture of an Exponential and an $\alpha-$Fr\'echet distributions.  Samples of size $N = 2^{17} = 131,072$ were
generated and a level $p=0.01$ and back--start parameter $b=4$ employed.
The left panel indicates the presence of a ``knee'' in the max--spectrum plot in one such mixture sample.
The automatic selection procedure identified well the location of the knee by selecting $j_1=9$ and the resulting
estimate $\what \alpha = 0.97$ is rather close to the nominal value of $\alpha=1$.  In the right panel,
we demonstrate the performance of the automatic selection procedure by using $1,000$ independent replications of the
mixture samples.  The histogram of the automatic choices for $j_1$ (left panel) indicates that most of the times
values close to the MSE--optimal one $j_1^{opt}=10$ were chosen.  The histogram of the resulting estimates of the
heavy--tail exponent (top--right graph in the left panel) is similar to the histogram corresponding to the MSE--optimal choice of $j_1$ (bottom--right in the left plot).  The slight bias in the histogram on the 
top--right is due to the fact that often slightly lower than the MSE--optimal values of $j_1$ were 
chosen by the automatic procedure.
More extensive analysis of this procedure is presented in Stoev et al.
 \cite{stoev:michailidis:hamidieh:taqqu:2006P}. 
 
\subsection{Data analysis}
 \label{s:data}

 We first discuss a popular insurance data set of $2,167$ fire losses in Denmark from 1980 to 1990.
 This data set has been studied extensively, see e.g.\ McNeil \cite{mcneil:1997}, Resnick \cite{resnick:1997d},  Lu and Peng \cite{lu:peng:2002} and Peng and Qi \cite{peng:qi:2004}.

 Figure \ref{fig:danish} displays the data, its corresponding Hill plot (bottom left) and its max--spectrum (bottom right).
 The max--spectrum yields an estimate $\what \alpha = 1.66$ obtained with an automatic selection of the scale
 $j_1$ by using a tunning parameter $p=0.01$ (see Section \ref{s:cut-off}), and the Hill plot yields an estimate
 $\what \alpha_H(k) = 1.39$ for $k = 1,000$.  This discrepancy  between the two methods is interesting since they
 yield comparable results in many typical models (see Section \ref{s:performance:typical}, above).  To explore
 further the significance of this difference, we resort to calculating confidence intervals.

 A particular advantage of the max--spectrum type estimators is that one can naturally obtain the
 following two types of confidence intervals for the parameters $H$ and $\alpha=1/H$: {\it (i)} based on the asymptotic normal distribution (see Proposition \ref{p:an-H-C})
 and {\it (ii)} based on a {\it permutation bootstrapping} procedure.  
 We will only briefly describe the procedure for obtaining permutation bootstrap confidence intervals. 
 Its theoretical analysis is outside the scope of the present paper.

\begin{figure}[h!]
\begin{center}
\includegraphics[width=5in]{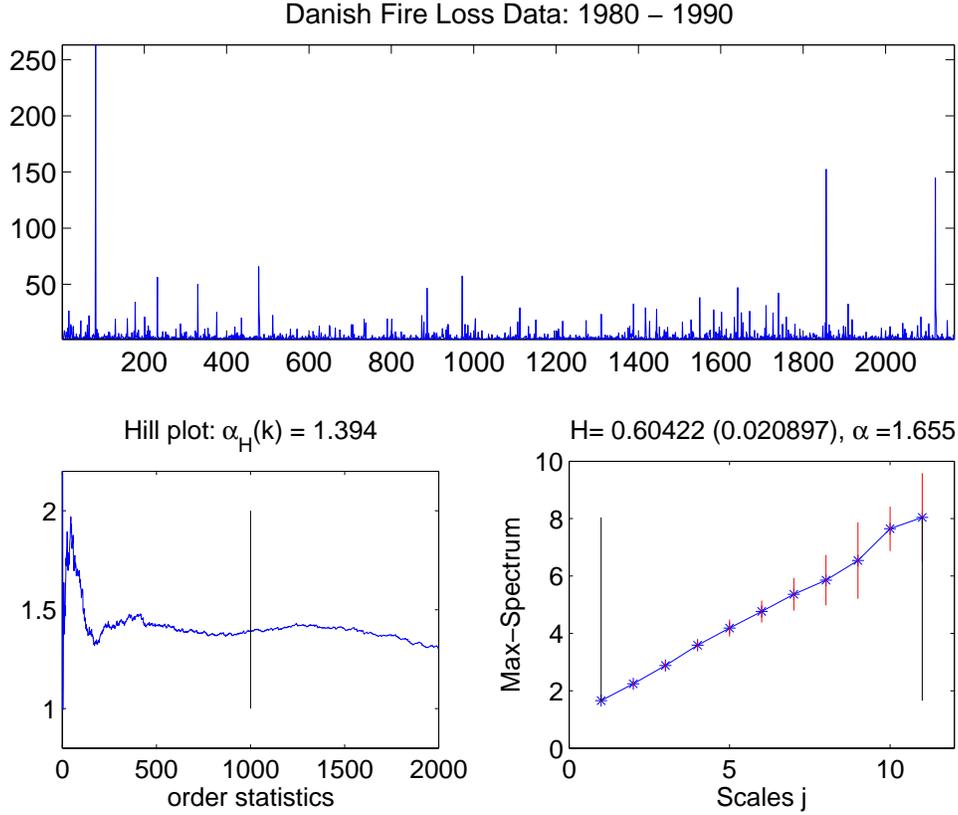} 
{\caption{\label{fig:danish} \small {\it Top panel:} time series of insurance losses due to 
  fire in Denmark from 1980 to 1990 losses (in million Danish krones).  Bottom left panel: the Hill plot
  of the fire loss data set.  {\it Bottom right:} the max--spectrum of the data.  Note that the Hill estimate is
  $\what \alpha_H(k) = 1.39,$ with $k=1,000$ and the max self--similarity estimate is $\what \alpha = 1.66$.
 }}
\end{center}
\end{figure}

 \medskip
{\bf Permutation bootstrap confidence intervals}   

\medskip 
 Given an i.i.d.\ sample $X(1),\ldots, X(N)$, generate $M$ independent random permutations
 $\pi_i:\{1,\ldots,N\} \to \{1,\ldots,N\}$, $i=1,\ldots,M$.  Then, construct the {\it permuted samples}
 $\wtilde X_i(1),\ldots,\wtilde X_i(N),\ i=1,\ldots,M$, where $\wtilde X_i(k) = X(\pi_i(k)),\ k = 1,\ldots,N$.
 Fix a range of scales $j_1<j_2\le \log_2 N$ and for each $i=1,\ldots,M$, compute the GLS max self--similarity estimator
 $\what H_i = \what H_i(j_1,j_2)$, from the permuted sample $\wtilde X_i(1),\ldots,\wtilde X_i(N)$.
 We will refer to the sample $\what H_i,\ i=1,\ldots, M$ as to the {\it permutation bootstrap} sample of the
 estimator $\what H = \what H(j_1,j_2)$, based on the original data set $X(1),\ldots, X(N)$.

 Observe that the statistics $\what H_i,\ i=1,\ldots,M$ are mutually dependent, since they are based on the
 original sample $X(1),\ldots, X(N)$.  However, since the $X(k)$'s are i.i.d.\ and the permutations $\pi_i$'s 
 are independent, we have that $\what H_i =^d \what H(j_1,j_2)$, for all $i=1,\ldots,M$.  One has moreover that 
 the sequence $\what H_i,\ i=1,\ldots, M$ is {\it exchangeable}.  This suggests using the permutation bootstrap
 sample $\what H_1,\ldots,\what H_M$ as a proxy to the sampling distribution of $\what H$.  We thus propose to use the empirical confidence interval based on the permutation bootstrap sample as a confidence interval
 for $H$.  Corresponding bootstrap confidence intervals for $\alpha = 1/H$ are obtained through the inversion
 method.

 \medskip
 Experience with several simulation experiments suggests the following conjecture.

 \begin{conjecture}  Let $\what H_i,\ i=1,\ldots,M$ be a permutation bootstrap sample of the estimator
 $\what H(j_1,j_2)$.  Consider the scales $j_1$, $j_2$ and the permutation sample size $M$ as functions of
 the sample size $N$, which tend to infinity as $N\to\infty$.  

 Under certain conditions on the rates of growth of $j_1,\ j_2$ and $M$, the empirical distribution of the 
 permutation bootstrap sample $\what H_i,\ i=1,\ldots,M$ yields asymptotically consistent confidence intervals for
 $\what H$.
 \end{conjecture}

\begin{figure}[t!]
\begin{center}
\includegraphics[height=2in,width=2.5in]{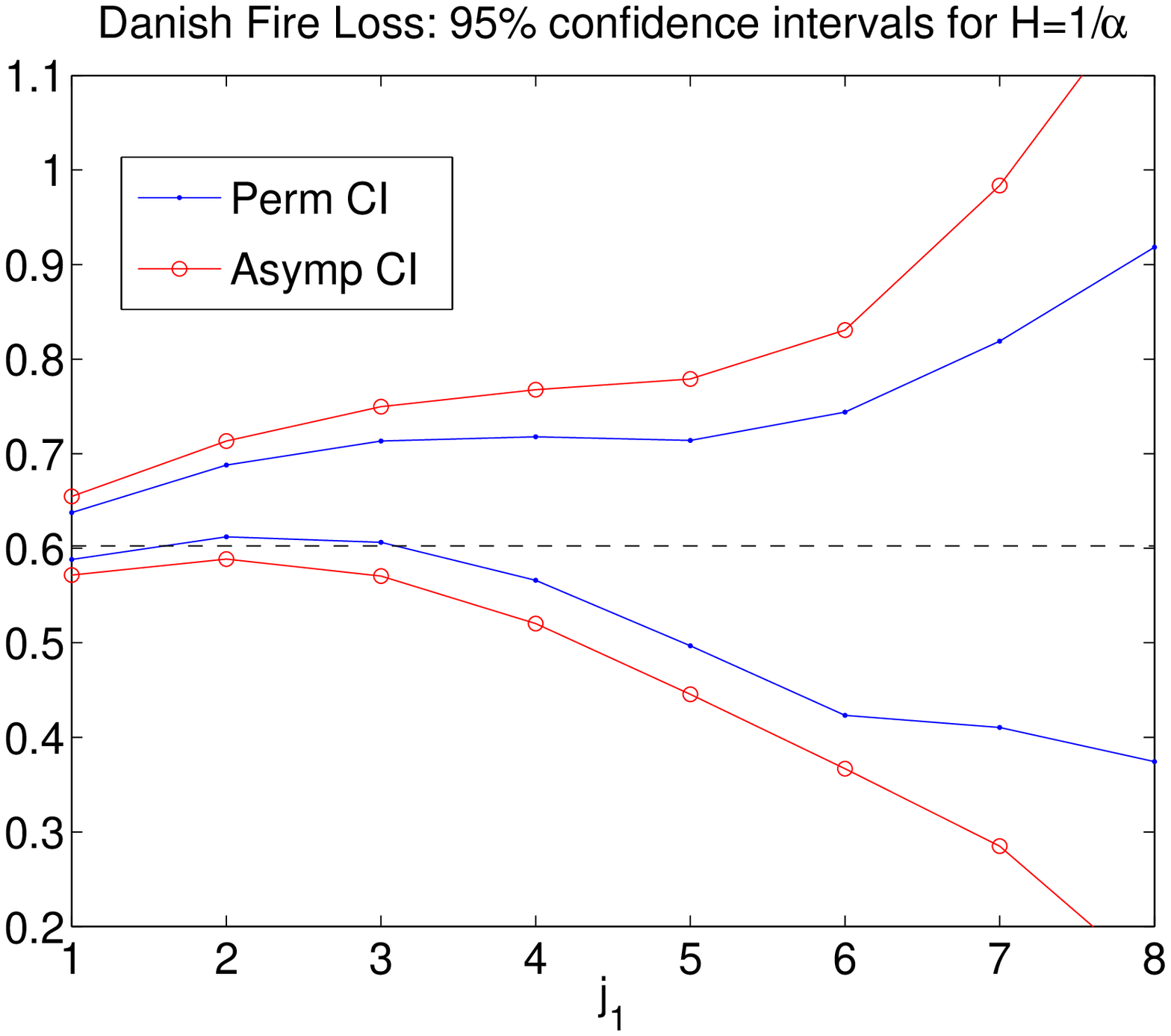} \hskip .2 in
\includegraphics[height=2in,width=2.5in]{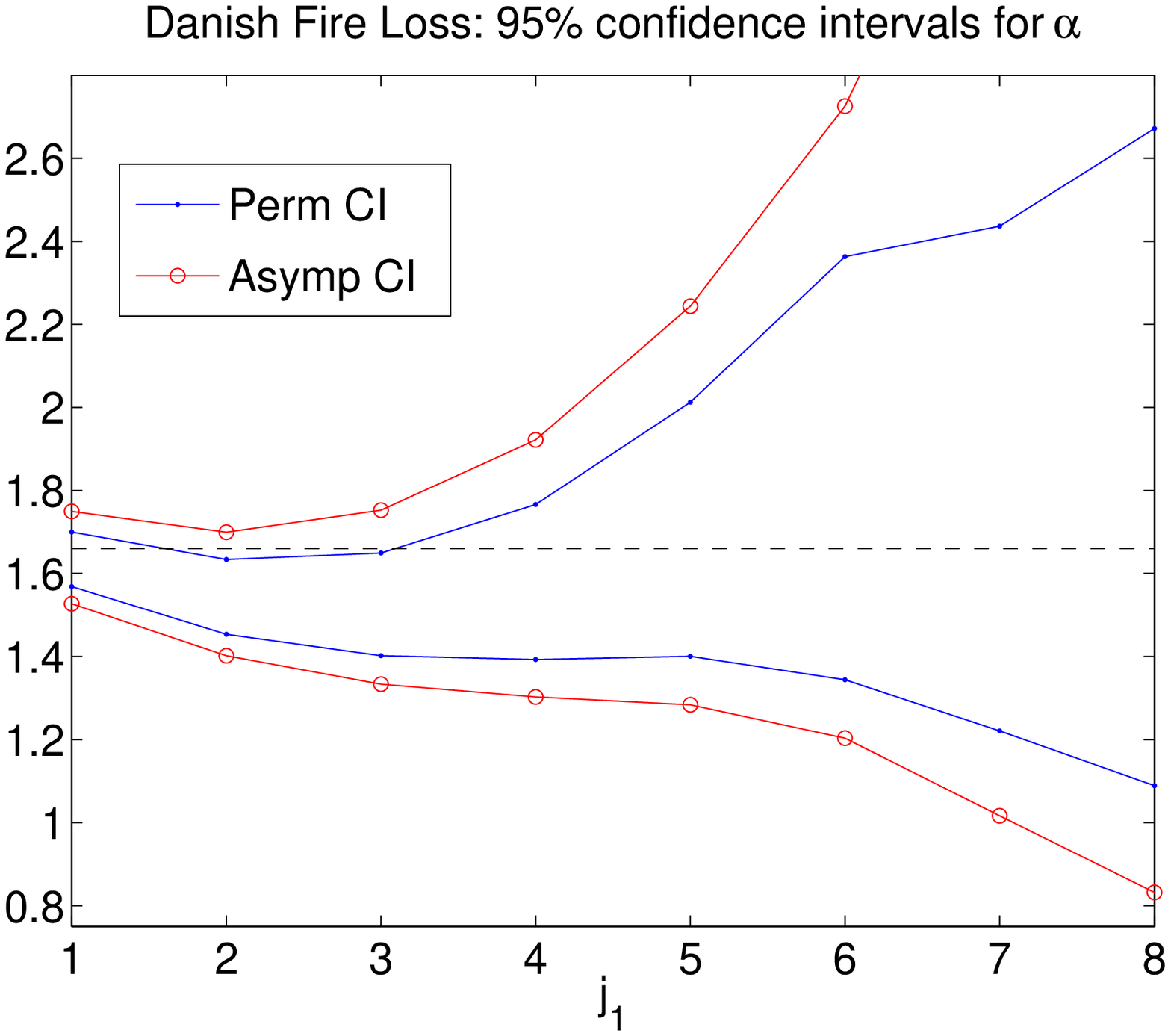}
{\caption{\label{fig:danish-ci} \small {\it Left panel:} $95\%$ confidence intervals
 for $H=1/\alpha$ based on: {\it (i)} permutation--bootstrap from $10,000$ independent
  permutations and {\it (ii)} asymptotic distribution for the max self--similarity estimators.
  {\it Right panel:} $95\%$ confidence intervals for $\alpha = 1/H$ obtained by inverting the 
  confidence intervals in the left panel.  The horizontal lines indicate the estimated value of
  $\what H = 0.6$ and $\what\alpha = 1/\what H = 1.66$ for $H$ and $\alpha$, respectively, obtained
   with the max self--similarity estimator in Figure \ref{fig:danish}.  }}
\end{center}
\end{figure}

\medskip
 Figure \ref{fig:danish-ci} displays $95\%$ confidence intervals for $H$ (left panel) and 
 $\alpha = 1/H$ (right panel) for the Danish fire loss data. Different scales $j_1$ were used and 
  $j_2$ was chosen as the largest available scale $11$.  The permutation confidence intervals (denoted by dots) 
  are obtained from
 $M=10,000$ random permutations and the asymptotic confidence intervals (denoted by circles) are obtained
 from the asymptotic variance in Proposition \ref{p:an-H-C} where the unknown value of $H$ was replaced by
 $\what H$.  To be able to compare the two types of intervals, we centered the asymptotic confidence intervals
 at the means of the permutation bootstrap samples $\what H_i,\ i=1,\ldots,M$.  Observe that although the two 
 procedures for constructing confidence intervals are different, they yield very similar results.  The permutation bootstrap
 intervals are always slightly more narrow than the asymptotic ones.  As Figure \ref{fig:danish} indicates,
 the use of scales $j_1 = 1$ and $j_2 = 11$ is acceptable.  The resulting permutation and asymptotic confidence intervals
 for $H$ are: $[0.5880, 0.6361]$ and $[0.5710, 0.6540]$, respectively.  They are consistent with, but considerably tighter than
 the likelihood--based intervals in Figure 8 of Lu and Peng \cite{lu:peng:2002} for the same data set.  This can be contributed 
 to the fact that the max--spectrum estimators and the Hill--type estimators are based on different 
 principles. 
 The performance of the permutation bootstrap and asymptotic confidence intervals is addressed in more
 detail in Stoev et al. \cite{stoev:michailidis:hamidieh:taqqu:2006P}.

\begin{figure}[h!]
\begin{center}
\includegraphics[width=5in]{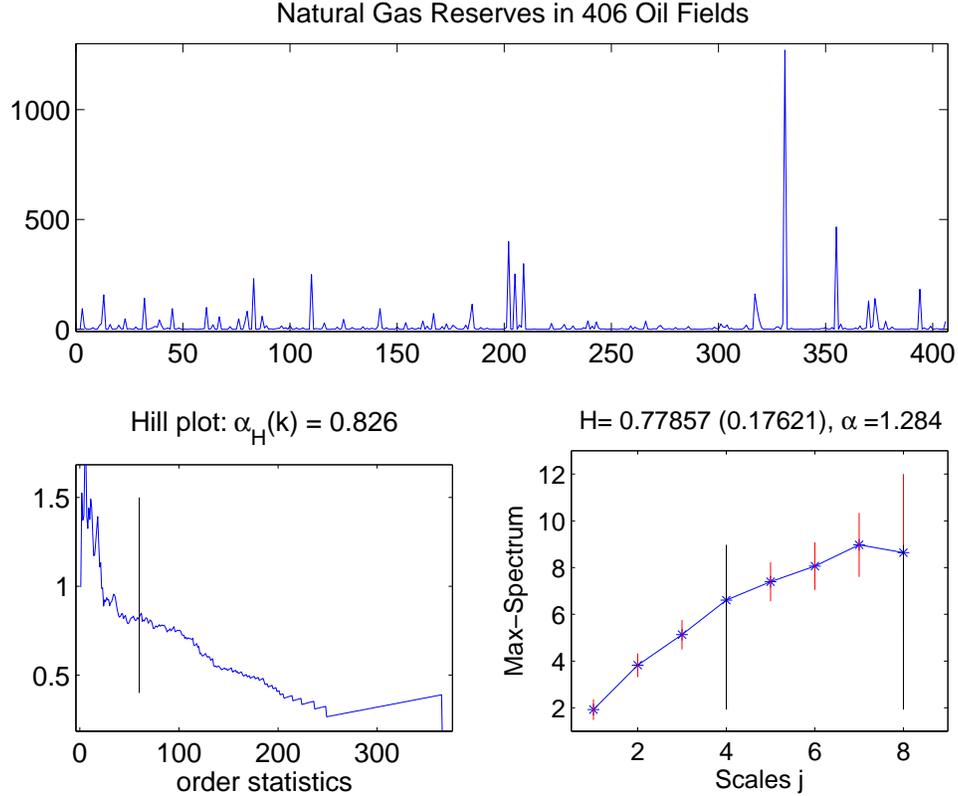} 
{\caption{\label{fig:gas} \small {\it Top panel:} randomly permuted sample of volumes natural gas reserves
 (in trillion cubic feet) found in $406$ provinces.  {\it Bottom left panel:} the Hill plot
  of the data set.  {\it Bottom right panel:} the max--spectrum of the data.  Note that the Hill estimate is
  $\what \alpha_H(k) = 0.826,$ with $k=60$ and the max self--similarity estimate is $\what \alpha = 1.284$.
 }}
\end{center}
\end{figure}

The second data set to be analyzed in this section consists of the volumes in trillion cubic feet of the $406$
largest natural gas world provinces. The data were obtained from Table 1 in
\cite{energy}. The study of the patterns in such data will help in  the
development of future natural gas resources leading to better
assessments of the reserve growth potential of the world's provinces. The
max self--similarity estimator, obtained from a typical randomly permuted sample
is $\what \alpha = 1.284$ (Figure \ref{fig:gas}).  Observe that the Hill plot shown 
in the bottom--left panel of Figure \ref{fig:gas} is very volatile and appears to stabilize
in a narrow range around $k=60$, where the resulting estimator is $\what\alpha_H(60) = 0.826$. 
Notice that the integer nature of the observations makes the Hill plot exhibit a saw-tooth like pattern
and hence difficult to obtain a good estimate for $\alpha$. Due to the discrepancy between the 
two methods, obtaining confidence intervals becomes particularly pertinent.

Permutation bootstrap and asymptotic confidence intervals for the max self--similarity estimators
for $H=1/\alpha$ and $\alpha$ are presented in Figure \ref{fig:gas-ci}.  As in Figure \ref{fig:gas-ci},
the asymptotic confidence intervals are slightly wider than the ones based on the permutation bootstrap. 
Observe that, contrary to the case of fire loss data in Figure \ref{fig:danish-ci}, the locations of the
confidence intervals for the gas data set stabilize only at scales $j\ge 4$.  This indicates that the
value $\what \alpha = 1.284$, obtained from the range of scales $j_1=4$ and $j_2$ in Figure \ref{fig:gas} is credible.
The fact that the resulting Hill estimate $\what \alpha_H(60) = 0.826$ is less than $1$ appears to be not
statistically significant, according to the confidence intervals in Figure \ref{fig:gas-ci}, which is in line
with the findings in de Sousa and Michailidis \cite{desousa:michailidis:2004}. 
This last fact and the volatility of the Hill plot suggest that the max self--similarity estimators can be
viewed as more reliable in this setting.

\begin{figure}[t!]
\begin{center}
\includegraphics[height=2in,width=2.5in]{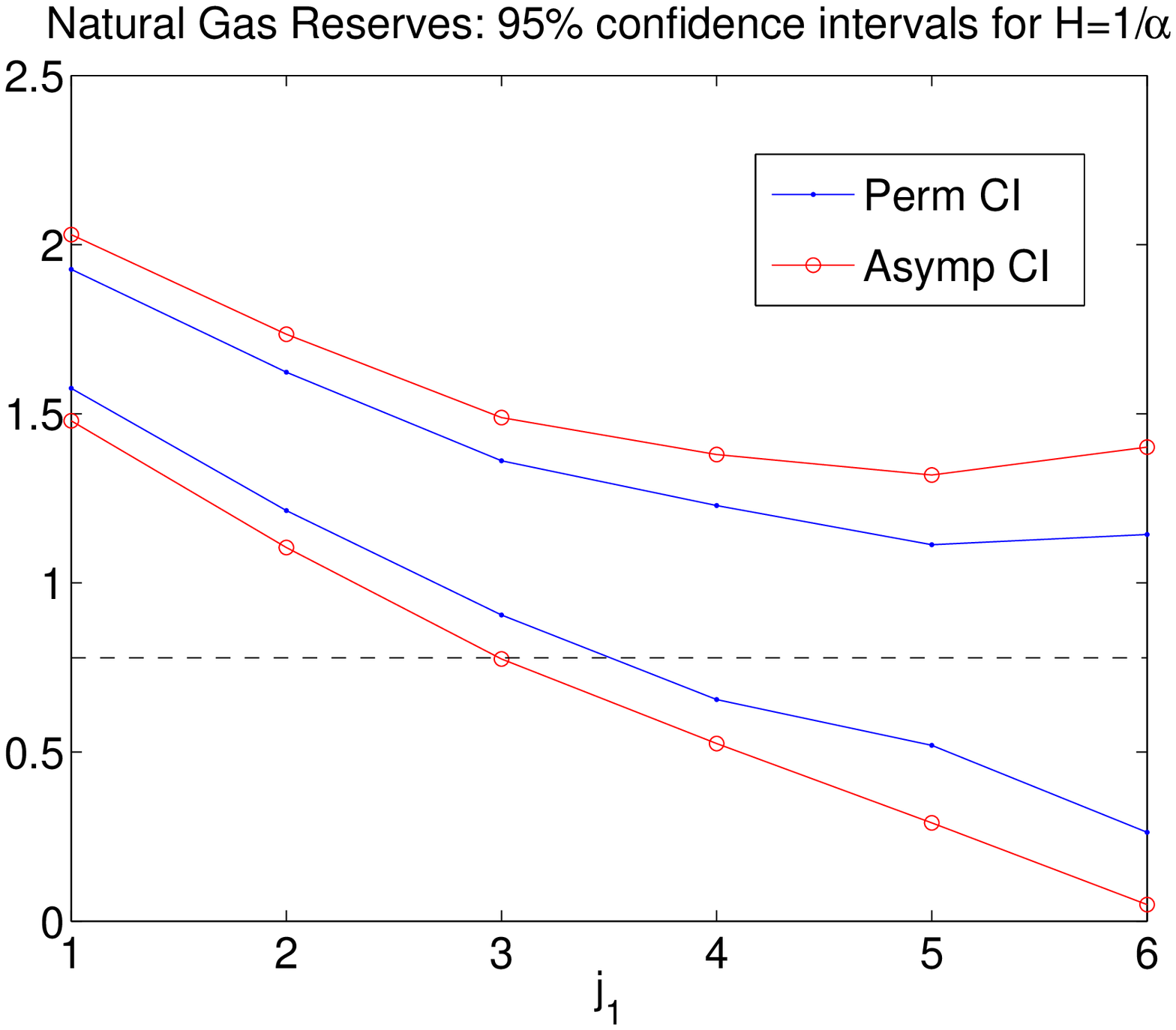} \hskip .2 in
\includegraphics[height=2in,width=2.5in]{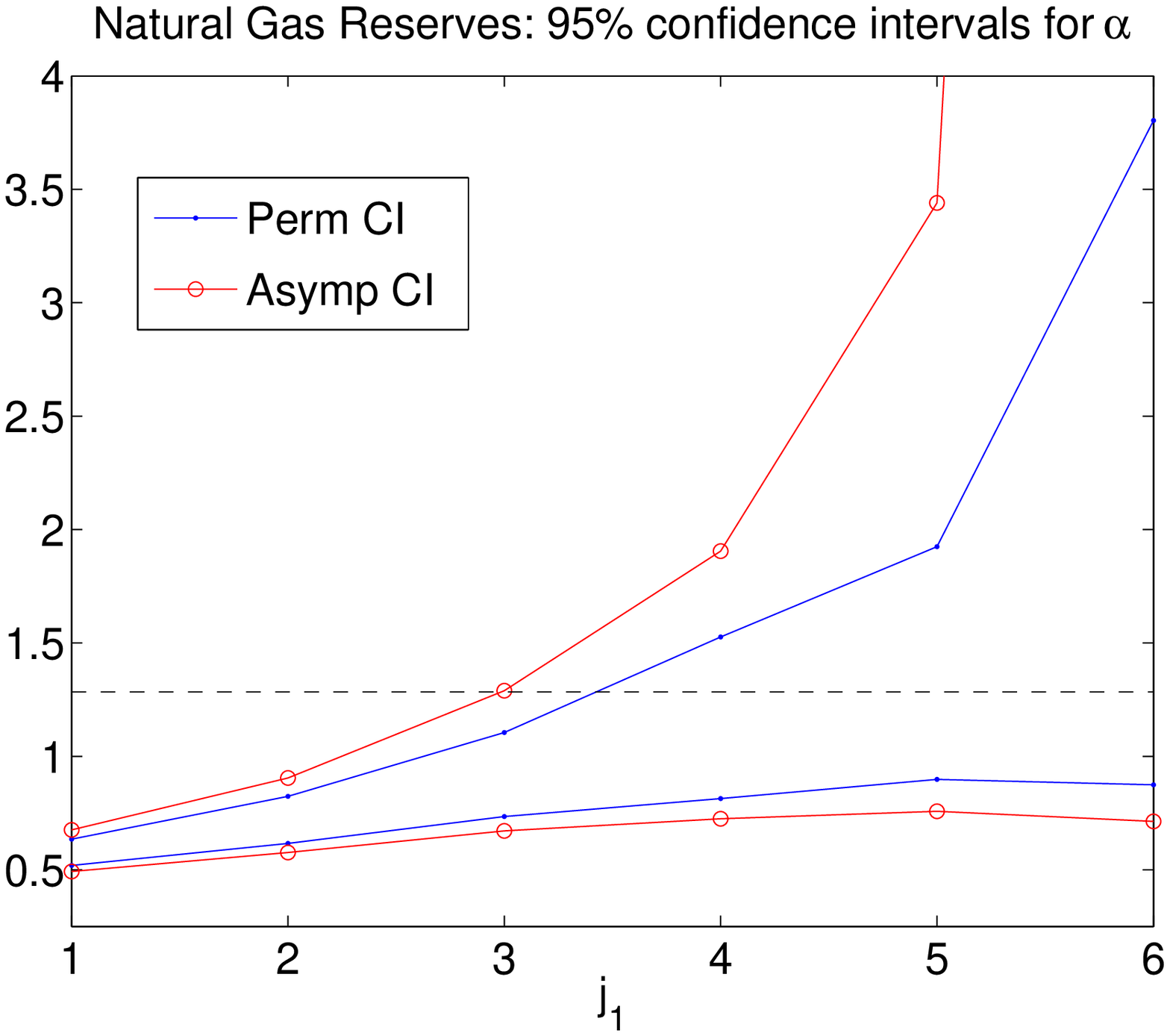}
{\caption{\label{fig:gas-ci} \small {\it Left panel:} $95\%$ confidence intervals
 for $H=1/\alpha$ based on: {\it (i)} permutation--bootstrap from $10,000$ independent
  permutations and {\it (ii)} asymptotic distribution for the max self--similarity estimators.
  {\it Right panel:} $95\%$ confidence intervals for $\alpha = 1/H$ obtained by inverting the 
  confidence intervals in the left panel.  The horizontal lines indicate the estimated value of
  $\what H = 0.78$ and $\what\alpha = 1/\what H = 1.28$ for $H$ and $\alpha$, respectively, obtained
   with the max self--similarity estimator in Figure \ref{fig:gas}.  }}
\end{center}
\end{figure}
\vskip 2in

\section{Concluding remarks}
 \label{s:conclusion}

In this paper, a new estimator for the tail exponent of a distribution 
was introduced and its asymptotic properties established. The estimator is based on block--maxima of the data and
can be visualized through a new graphical device called the max--spectrum plot. Numerical work shows that
compared to the widely used Hill estimator, the max self--similarity estimator performs
competitively in the case of the Pareto distribution and it outperforms the
Hill estimators in the cases of the stable, Fr\'echet and certain t-distributions. In practice, the max--spectrum
plot is less volatile than the classical Hill plot.  Thus, the max self--similarity estimator can be used
in situations where the Hill plot fails or when it is hard to interpret.  Finally, the fact that the estimator
is based on block maxima makes it particularly suitable for time series data, a topic discussed in a companion
paper Stoev et al.\ \cite{stoev:michailidis:hamidieh:taqqu:2006P}.

\medskip
\noindent
{\bf Acknowledgments:}  The authors would like to thank Professor Robert Keener
for suggesting the proof of Lemma \ref{l:int-by-parts} and for many useful
suggestions.  We also thank Kamal Hamidieh for stimulating discussions on the automatic selection
procedure of scales. 
The work was supported in part by a grant from the Horace H.\ Rackham School of Graduate
Studies at the University of Michigan (SS), and by NSF grants CCR-0325571,
DMS-0505535 (GM) and DMS-0505747 (MT).

\section{Appendix: auxiliary results and tables}

\subsection{Auxiliary results} 

 We briefly review some properties of the $\alpha-$Fr\'echet distributions used above.

 \begin{definition}\label{d:Frechet} A random variable $Z$ is said to have an $\alpha-$Fr\'echet distribution,
 if
 \beq\label{e:Z-x}
  \P\{Z \le x \} = \left\{\begin{array}{ll}
                            \exp\{ -\sigma^\alpha x^{-\alpha}\} &,\ x>0 \\
                             0 &,\ x\le 0,
                          \end{array}
                   \right.
 \eeq
 with $\sigma>0$.  The parameter $\sigma$ is referred to as the scale coefficient of $Z$.  The random
 variable $Z$ is said to be standard $\alpha-$Fr\'echet if $\sigma=1$.
 \end{definition}

\medskip Let $Z$ be an $\alpha-$Fr\'echet variable with scale coefficient $\sigma>0$.
 The next properties follow directly from Relation \refeq{Z-x}.  

 \medskip
\noi {\bf Properties}
\begin{enumerate}
 \item {\it (scale family)}  For all $c>0$, the random variable $cZ$ is $\alpha-$Fr\'echet and
 has scale coefficient $c\sigma$.

 \item {\it (heavy tails)} The Taylor expansion of the exponential around the origin implies that
 \beq\label{e:Z-tail}
  \P\{Z > x\} = 1 - e^{-\sigma^\alpha x^{-\alpha}} \sim \sigma^\alpha x^{-\alpha},\ \ \mbox{ as } x\to\infty.
 \eeq
 
 \item {\it (moments)} In view of \refeq{Z-tail}, for all $p>0$,
 $$
 \E Z^p <\infty \ \ \mbox{ if and only if }\ \ p<\alpha.
 $$
One has moreover, that $\E Z^p = \sigma^p \Gamma(1-p/\alpha),\ p\in(0,\alpha)$, with
$\Gamma(x)=\int_0^\infty u^{x-1}e^{-u}du,\ x>0$.

 \item {\it (log--moments)} For all $p>0$, the moments $\E |\ln Z|^p$ are finite.  This follows from the 
 fact that $\xi:= \alpha\ln(Z/\sigma)$ has the Gumbel distribution, i.e.\ $\P\{\xi\le x\} = \exp\{ -e^{-x}\},\ x\in\bbR$.
 See also Corollary \ref{c:logs} below.

 \item {\it (power transformations)} For any $p>0$, the random variable $Z^p$ is $\alpha/p-$Fr\'echet with 
 scale coefficient $\sigma^p$.  Consequently, if $Z_1$ is a standard $1-$Fr\'echet variable, then 
 $$
  Z := Z_1^{1/\alpha} 
 $$
 is standard $\alpha-$Fr\'echet, for all $\alpha>0$.
\end{enumerate}

\medskip
\noi The $\alpha-$Fr\'echet distributions are also {\it max--stable} in the following sense.

\begin{definition}\label{d:max-stable}
 A random variable $Z$ is said to be max--stable, if for all $a, b>0$ there exist $c>0,\ d\in\bbR$, such that
$$
\max\{ a Z', b Z''\} \stackrel{d}{=} c Z + d,
$$
where $Z'$ and $Z''$ are independent copies of $Z$ and where $=^d$ means equality in distribution.
\end{definition}

In particular, by \refeq{Z-x}, one gets that if $Z(1),\ldots, Z(n),\ n\in\bbN$ are i.i.d.\  $\alpha-$Fr\'echet, then
\beq\label{e:Z-max-stable}
  Z(1)\vee \cdots \vee Z(n) \stackrel{d}{=} n^{1/\alpha} Z(1). 
\eeq
This last relation shows that a sequence of i.i.d.\  $\alpha-$Fr\'echet variables is also max self--similar with
parameter $H=1/\alpha$ (see Definition \ref{d:max-ss} above).  Relation \refeq{Z-max-stable} served
as the main motivation to define the max self--similarity estimators in Section \ref{s:max-ss} above.

The class of max--stable distributions in the sense of Definition \ref{d:max-stable} above includes, in addition
to the Fr\'echet, only the classes of negative Fr\'echet and the Gumbel laws.  These three classes of distributions
are the only distributions arising in the limit of maxima of i.i.d.\  variables under appropriate normalization
(see e.g.\ Proposition 0.3 in Resnick \cite{resnick:1987} and also
 Leadbetter, Lindgren and Rootz\'en \cite{leadbetter:lindgren:rootzen:1983}).

\medskip
\noi The following integration by parts formula is used in the proof of Theorem \ref{t:rate}.

\begin{lemma}\label{l:int-by-parts}
 Let $f:[a,b]\to \bbR,\ a,b\in\bbR $ be an absolutely continuous function, that is, 
 $f(x) = f(a) + \int_a^x f'(u) du,$ for some Lebesgue integrable $f'(x),\ x\in[a,b]$.  
 Then, for any c.d.f.\ $G(x)$, we have 
\beq\label{e:l:int-by-parts}
 \int_a^{b} f(x) d G(x) = f(b)G(b) - f(a) G(a) - \int_a^b G(x) f'(x) dx.
\eeq
\end{lemma}
\begin{proof} Since $f(x) = f(a) + \int_a^b f'(u) 1_{[a,x)}(u)du$, we have that
$$
 \int_a^{b} f(x) d G(x) = f(a)G(b) - f(a)G(a) + \int_a^{b} {\Big(}\int_a^{b} f'(u) 1_{[a,x)}(u)du {\Big)} d G(x).
$$
An application of Fubini's theorem yields
\begin{eqnarray*}
& & f(a)G(b) - f(a)G(a) + \int_a^b f'(u) (G(b) - G(u)) du \\
 & &  \ \ \ \ \ \ \ \  =  f(a)G(b) - f(a)G(a) + (f(b) - f(a))G(b) - \int_a^b f'(u) G(u) du.
\end{eqnarray*}
Observe that the right--hand sides of the last expression and Relation \refeq{l:int-by-parts} coincide.
 $\Box$
\end{proof}

\subsection{Tables}

\include{stable_hill_table_1}

\include{t_dist_hill_table_1}

 \begin{table}[h!]
 \begin{center}
  \begin{tabular}{||c||c|c|c|c|c||}
 \hline
  $\psi(i)$ &  $i+0$ & $i+1$ & $i+2$ & $i+3$ & $i+4$ \\
 \hline\hline
 $i=0$  &  $3.423696$ & $2.211864$ & $1.387207$ & $0.846734$ & $0.504666$ \\
 \hline
 $i=5$  &  $0.294581 $ & $0.168963 $ & $0.095563$ & $0.053288 $ & $0.029470$ \\
 \hline
 $i=10$ &  $0.016072$ & $ 0.008755 $ & $ 0.004756 $ & $ 0.002552$ & $ 0.001405 $ \\
 \hline
 $i=15$ &  $ 0.000709 $ & $0.000335$ & $ 0.000175 $ & $0.000097 $ & $0.000032$ \\
 \hline\hline
 \end{tabular}
 \end{center}

 \caption{\label{tab:psi} {\small We present here numerical approximations of the values
 $\psi(i)$, $i=0,1,\ldots,19$ involved in the expression of the covariance matrices 
 $\Sigma_\alpha(j_1,j_2;N)$ in \refeq{p:cov} (see also \refeq{psi}).  We used Monte Carlo
 simulations with $10, 000, 000$ independent pairs of $1-$Fr\'echet variables.  To reduce the 
 variance of the estimates we used ``bagging''.  That is, the Monte Carlo simulations were repeated
 independently $1, 000$ times and then the resulting means were taken as the final estimates 
 reported in the table above.}}
 \end{table}

 \begin{table}[h!]
\begin{center}
  \begin{tabular}{||c||c|c|c|c|c|c|c|c|c|c||}
 \hline 
 $j$ & $ 2 $ & $ 3 $ & $ 4 $ & $ 5 $ & $ 6 $ & $ 7 $ & $ 8 $ & $ 9 $ & $ 10 $ & $ 11$ \\
\hline\hline
$\sqrt{c_w(j)}$  & $1.417$ & $ 0.802$ & $ 0.515$ & $ 0.346$ & $ 0.238$ & $ 0.166$ & $ 0.116$ & $ 0.082$ & $ 0.058$ & $ 0.041$\\
\hline
$j+=10$  & $0.029$ & $ 0.020$ & $ 0.014$ & $ 0.010$ & $ 0.007$ & $ 0.005$ & $ 0.004$ & $ 0.003$ & $ 0.002 $ & $ 0.001$\\
\hline\hline
$\sqrt{2^j c_w(j)}$  & $2.834$ & $ 2.267 $ & $ 2.060  $ & $ 1.960$ & $1.905 $ & $ 1.875 $ & $1.857  $ & $ 1.847 $ & $ 1.841 $ & $ 1.837 $\\
\hline
$j+=10$  & $1.835 $ & $1.834 $ & $ 1.834 $ & $ 1.833$ & $1.833 $ & $1.833 $ & $ 1.833 $ & $1.833 $ & $1.833 $ & $ 1.833 $\\
\hline\hline
 \end{tabular}
 \end{center}

 \caption{\label{tab:cw} {\small 
We present here numerical estimates of the constants $c_w$ involved in the asymptotic 
 variances in Proposition \ref{p:an-H-C} above.  Here, we use $j_1=1$, for simplicity, and display $20$ different values
 corresponding to $j_2 = j = 2,\ldots,21$. For convenience, we present $\sqrt{c_w}$ together with $\sqrt{2^{j_2}c_w}$ where the
 latter constant is useful if one normalizes in \refeq{an-H-C} by using $\sqrt{N_r}$ instead of $\sqrt{N_{j_2+r}}$.}}
 \end{table}

\section* 
 \small
 
\bibliography{/afs/umich.edu/user/s/s/sstoev/doc/articles/mst-bibfile.bib}

\end{document}